\documentclass[journal]{new-aiaa}
\usepackage[utf8]{inputenc}
\usepackage{textcomp}

\usepackage{graphicx}
\usepackage{amsmath}
\usepackage[version=4]{mhchem}
\usepackage{siunitx}
\usepackage{longtable, tabularx, booktabs, multicol, multirow}
\usepackage{comment}
\usepackage{arydshln}
\usepackage{float}
\usepackage{algorithm}
\usepackage{subcaption}
\usepackage{algorithmic}
\usepackage[normalem]{ulem}
\usepackage{tikz}
\usetikzlibrary{
  arrows.meta,
  positioning,
  shapes.geometric,
  calc
}

\newcommand{\e}[1]{\times 10^{#1}}  

\title{Convex Trajectory Optimization via Monomial Coordinates Transcription for Cislunar Rendezvous}

\author{O. Regantini\footnote{Ph.D. Student, Department of Aerospace Science and Technology, Via La Masa 34, \href{mailto:omar.regantini@polimi.it}{omar.regantini@polimi.it}}, E. R. Burnett\footnote{Marie Skłodowska-Curie Postdoctoral Fellow, Department of Aerospace Science and Technology, Via La Masa 34, \href{mailto:ethanryan.burnett@polimi.it}{ethanryan.burnett@polimi.it}, Member AIAA}, A. Rizza\footnote{Ph.D., Department of Aerospace Science and Technology, Via La Masa 34, \href{mailto:antonio.rizza@polimi.it}{antonio.rizza@polimi.it} Member AIAA}, A. Morselli\footnote{Assistant Professor, Department of Aerospace Science and Technology, Via La Masa 34,\href{mailto:alessandro.morselli@polimi.it}{alessandro.morselli@polimi.it}, Member AIAA} and F. Topputo\footnote{Professor, Department of Aerospace Science and Technology, Via La Masa 34,\href{mailto:francesco.topputo@polimi.it}{francesco.topputo@polimi.it}, Senior Member AIAA}}
\affil{Polytechnic University of Milan, 20156 Milan, Italy}

\begin{document}

\maketitle

\begin{abstract}
This paper proposes a nonlinear guidance algorithm for fuel-optimal impulsive trajectories for rendezvous operations close to a reference orbit. The approach involves overparameterized monomial coordinates and a high-order approximation of the dynamic flow precomputed using differential algebra, which eliminates the need for real-time integration. To address non-convexity in the monomial coordinate formulation of the guidance problem, sequential convex programming is applied. Using the methodology presented in this paper, repeatedly evaluating the nonlinear dynamics is not necessary, as in shooting or collocation methods. Instead, only the monomial equations require updating between iterations, drastically reducing computational burden. The proposed algorithm is tested in the circular restricted three-body problem framework with the target spacecraft on a near-rectilinear halo orbit. The results demonstrate stability, efficiency, and low computational demand while achieving minimal terminal guidance errors. Compared to linear methods, this nonlinear convex approach exhibits superior performance in open-loop propagation of impulsive maneuvers in cislunar space, particularly in terms of accuracy. These advantages make the algorithm an attractive candidate for autonomous onboard guidance for rendezvous operations in the cislunar domain.
\end{abstract}

\section{Introduction}
Orbital rendezvous is a critical technology in modern spaceflight. This capability is key in advanced space operations, such as crew transfers, cargo exchange, refueling, and satellite repairs. {These missions consist of a sequence of maneuvers and a controlled trajectory that brings the active spacecraft (\emph{chaser}) close to the passive vehicle (\emph{target}).} In \cite{woffinden_navigating_2007}, a comprehensive review of the history, technological development, and evolution of rendezvous missions is given. As future missions extend farther from Earth, minimizing the effects of communication latency and potential data loss becomes essential. Hence, autonomous Guidance, Navigation, and Control (GNC) systems become necessary for several new challenging tasks with minimal ground intervention. Given the computational and hardware limitations of onboard computers, onboard autonomous GNC systems must balance computational and power efficiency with robust performance. 


The field of GNC has undergone a fundamental shift, with a growing emphasis on numerical and computational methods for onboard execution \cite{lu_introducing_2017}. This trend has transformed how intricate guidance problems are approached, moving away from traditional control laws and focusing on optimization methods that accommodate nonlinear dynamics and constraints. In space guidance, optimal control problems (OCPs) have historically been solved using either indirect or direct methods \cite{betts_survey_1998}. Indirect methods are based on the necessary conditions derived from the calculus of variation and require explicit formulation and integration of the associated costate dynamics \cite{berkovitz_optimal_1974, liberzon_calculus_2012}. Primer vector theory \cite{primer-vector-1963} has been widely used in this framework to determine fuel-optimal trajectories for both impulsive and continuous thrust. 
On the other hand, direct methods transcribe the continuous-time OCP into a nonlinear programming problem through discretization of state and control variables\cite{directOptmization}. Among these, direct transcription and collocation (Hermite-Simpson, high-order Gauss-Lobatto, and pseudospectral method) have been widely adopted and studied in the optimal control problem of space trajectories \cite{survey-topputo}. Although indirect methods offer high accuracy, they require complex analytical derivations and careful selection of initial guesses of costate variables, whereas direct methods, though more straightforward, face convergence challenges and are computationally expensive. Recently, convex optimization has emerged as a particularly powerful framework in this context \cite{mao_successive_2017, mao_tutorial_2018}, overcoming the limitations of canonical methods by providing a reliable, computationally efficient approach that works well in real-time autonomous scenarios even with complex dynamics and constraints. Indeed, convex optimization is a robust approach that ensures global convergence in polynomial time \cite{boyd_convex_2004, nesterov_lectures_2018}. However, many optimal control problems are inherently non-convex due to the complex nonlinear dynamics and constraints \cite{acikmese_convex_2007, wang_minimum-fuel_2018}. Consequently, in the case of non-convex problems, various methods have been proposed to turn a non-convex optimization problem into a convex one, such as lossless convexification and {sequential convex programming} (SCP) \cite{liu_survey_2017}. Collocation-based SCP formulations \cite{hofmann_performance_2021}, both for continuous-thrust and impulsive control, enforced the system dynamics algebraically at the collocation points, thus avoiding real-time numerical integration of the equation of motion. In contrast, alternative approaches, such as first-order hold or state transition matrix (STM)-based impulsive SCP \cite{benedikter2019convexoptimizationlinearimpulsive, impulsiveSTMSCP_rizza}, rely on the integration of the equation of motion, along with the variational equations, at each SCP iteration. While STM-based methods incur higher computational cost, they generally provide increased accuracy in modelling the system response compared to collocation-only formulations \cite{hofmann_performance_2023}.

Several recent efforts to address the limited on-board computational capability have focused on using machine learning to solve the guidance problem. 
During training, the system learns how to map spacecraft states to the corresponding controls. These initial efforts are encouraging, but the training process for these methods can be time-consuming, highly sensitive to parameter tuning, and has no performance guarantees. An alternative research direction exploits high-order expansion of the dynamics to enable the approximation of the trajectory incorporating the system nonlinearities while maintaining computational tractability. For example, in \cite{boone_orbital_2021, boone_optimal_2021}, a predictor-corrector guidance strategy is presented for station-keeping in an unstable halo orbit within the lunar domain, leveraging state transition tensors (STTs) \cite{park_nonlinear_2006}. These STTs can be thought of as higher-order extensions of the STM. Although the integration through variational equations of these higher-order terms can be computationally demanding, once computed for a reference trajectory, the spacecraft's state can be approximated by simple algebraic operations, eliminating the need for onboard integrations. In \cite{burnett_rapid_2025}, a high-order analytical expansion of the relative dynamics is combined with SCP, yielding a rapid and stable approach for on-board implementation with application to Earth's orbit. Though this approach is open-loop, it remains computationally efficient and can solve for a new optimal guidance solution in a closed-loop fashion as needed, based on the current state estimate. 

This paper introduces a nonlinear convex method for onboard rendezvous guidance that avoids real-time integration, instead, it requires few nonlinear operations to be executed on board, it is faster, and it offers greater assurance of accuracy than classical transcription methods. By employing a high-order expansion of the dynamics, the method avoids linearization or model approximations, which can introduce significant errors in the case of large spacecraft separations and long timescales, especially for highly nonlinear environments like the cislunar domain. These high-order expansions are computed using differential algebra (DA) \cite{berz_chapter_1999} methods, which automatically compute high-order derivatives. The optimization problem is efficiently formulated using a specialized set of overparameterized monomials based on the initial conditions of the problem, resulting in a novel transcription scheme \cite{burnett_rapid_2025}. Although the work in \cite{burnett_rapid_2025} deals with a fixed-final time rendezvous scenario for Earth orbit, this work introduces a free-final time formulation and applies their approach to a more complex dynamical domain. 
Within the context of the circular restricted three-body problem (CR3BP), several optimal control approaches have been proposed to solve impulsive rendezvous problems. Notably, primer vector-based methods have been used to compute orbit-to-orbit impulsive transfers \cite{bucchinio-primer-vector} as well as rendezvous trajectories by linearizing the dynamics about libration points \cite{serra-primer-vector}. Complementing these approaches, the method presented in this paper uses an SCP optimization framework, enabling robust and efficient guidance computation without relying on local linear models or real-time integration.
This approach is thoroughly tested in cislunar space, with the target spacecraft placed into a {near rectilinear halo orbit} (NRHO). NRHOs are particularly relevant in the context of the lunar Gateway mission, as they offer long-term stability and low-energy access to the lunar surface, making them key enablers for sustained lunar exploration and logistics \cite{lee_white_2019}. However, rendezvous operations in NRHO are significantly more challenging due to the highly nonlinear and sensitive dynamics of the CR3BP framework.  The algorithm is evaluated across various test cases, each characterized by different levels of nonlinearity and varying initial target-chaser separations. Such scenarios offer a wide range of testbeds to analyze the robustness and performance of algorithms in various mission conditions operating in complex space environments.
In addition to numerical evaluation, the optimality of the computed solutions is verified \textit{a posteriori} using primer vector theory, and both computational efficiency and solution quality are quantitatively compared against a classical SCP formulation.


This paper is organized as follows: Section \ref{sec:mon para} covers the basic theory of monomial parameterization and the tool to rapidly and automatically extract high-order expansions.  Section \ref{sec:dyn model} outlines the assumption and the model used to describe the system's dynamics, and once the dynamics is established, the focus shifts to formulating the optimal control problem using the monomial matrix in Section \ref{sec:ocp}. Finally, in Section \ref{sec:results}, the performance of the proposed nonlinear convex method will be evaluated by applying it to various test cases within the CR3BP, and concluding remarks are given in Section \ref{sec: conclusion}.

\section{Theoretical Background} \label{sec:mon para}
\subsection{Monomial Parameterization}
Consider the truncated $m$-order Taylor Series Expansion (TSE) of a $n$-variables scalar function $f:\mathbb{R}^n \rightarrow \mathbb{R}$, where $f\in C^m(\mathbb{R}^n)$,  expanded around the point $\boldsymbol{\bar{x}}$:
\begin{equation} \label{eq :TSE def}
    f(\boldsymbol{x}) = f(\bar{\boldsymbol{x}}) + \sum_{\gamma_1=1}^n\dfrac{\partial f}{\partial x_{\gamma_1}}\bigg|_{\bar{\boldsymbol{x}}}\delta x_{\gamma_1} + \dfrac{1}{2!}\sum_{\gamma_1=1}^n\sum_{\gamma_2=1}^n\dfrac{\partial^2 f}{\partial x_{\gamma_1} \partial x_{\gamma_2}}\bigg|_{\bar{\boldsymbol{x}}}\delta x_{\gamma_1} \delta x_{\gamma_2} + \dots 
\end{equation}
This series can be compactly reformulated using matrix notation as:
\begin{equation} \label{eq:mon_para_def}
    f(\boldsymbol{x}) \approx f(\boldsymbol{\bar{x}}) + \Psi_m \boldsymbol{c}_m
\end{equation}
Here,  $\Psi_m$ is the row vector collecting all evaluated partial derivatives of the expansion, and $\boldsymbol{c}_m$ gives the unique monomial combinations of the components of $\delta \boldsymbol{x}$ up to order $m$. The terms in $\Psi_m$ and $\boldsymbol{c}_m$ are ordered similarly to \citet{giorgilli_methods_2011}. 
This concept extends naturally to vector functions, such as the flow of a dynamical system, where $\Psi_m$  becomes a matrix with $n$ rows, each row representing the Taylor series expansion of the corresponding component of the vector function \cite{burnett_rapid_2025}. In this framework, each monomial $\boldsymbol{c}_m$ corresponds to a basis function to describe the system's flow that is usually numerically encoded in the STTs. Because the two approaches are equivalent, we use the terms \textit{monomial parametrization}, \textit{state transition tensor}, \textit{Taylor expansion}, or \textit{map} interchangeably to refer to the matrix $\Psi_m$. The reader is encouraged to consult Appendix \ref{app:stt and mon matrix equi} for an in-depth analysis of this equivalency. In the context of a dynamical system, $\Psi_m$ is time-dependent, meaning that the matrix $\Psi_m(t_0, t_i)$, which approximates the state at time $t_i$ from an initial deviation at time $t_0$, differs from the monomial matrix $\Psi_m(t_0, t_j)$ used to approximate the state at another time $t_j \neq t_i$. Therefore, while for an algebraic function, $\Psi_m$ is unique, in the case of a dynamical system flow, $\Psi_m(t_0, t_i)$ varies over time.

The vector $\boldsymbol{c}_m$ contains $K_m$ elements, where $K_m$ represents the number of unique combinations of differentials up to order $m$:
\begin{equation}
    K_m = \dfrac{(n+m)!}{n!m!} - 1
\end{equation}
Among these $K_m$ elements, the first $n$ are independent and form the linear part, denoted as $\boldsymbol{c}_1$, while the remaining $(K_m - n)$ elements are {defined} by this linear part. These constraints define a nonlinear relationship that restricts $\boldsymbol{c}_m$ to lie on an $n$-dimensional embedded submanifold $\mathcal{C}^{(n, m)}$ of $\mathbb{R}^{K_m}$. As a result, the entire vector $\boldsymbol{c}_m$ can be uniquely determined from $\boldsymbol{c}_1$ through an analytic mapping $\boldsymbol{E}_m:\mathbb{R}^n\rightarrow\mathbb{R}^{K_m}$ \cite{burnett_rapid_2025}:
\begin{equation}
\label{eq:mon_pro}
    \boldsymbol{c}_m = \boldsymbol{E}_m(\boldsymbol{c}_1)
\end{equation}
The function $\boldsymbol{E}_m$ defines a generic $j$-th element of $\boldsymbol{c}_m$, $c_{j, m},$ as:
\begin{equation}
    c_{j, m} = \prod_{k=1}^n c_{k, 1}^{\sigma_{jk} }
\end{equation}
where $c_{k, 1}$ denotes the $k$-th component of the linear part, and $\sigma_{jk}$ is element in the position $(j, k)$\footnote{In this work, $m$ denotes the expansion order, while $j$ represents a generic index. In contrast, in \cite{burnett_rapid_2025}, $j$ indicates the order of the expansion.}  of  $\Sigma$, a $K_m \times n$ matrix which collects the exponents of the monomials. The first $n$ rows of  $\Sigma$ are related to the linear part of the expansion, and thus the top $n\times n$ sub-matrix of $\Sigma$ is the identity. {An} algorithm for defining $\Sigma$ is detailed in Appendix A of \cite{burnett_rapid_2025}. 

Given a point $\boldsymbol{c}_m \in \mathcal{C}^{(n, m)}$, any infinitesimal variations $d\boldsymbol{c}_m$ must belong to its tangent space, and an infinitesimal deviation that lies in the tangent plane satisfies the following relation:
\begin{equation}
    d\boldsymbol{c}_m = \dfrac{\partial \boldsymbol{c}_m}{\partial \boldsymbol{c}_1}d\boldsymbol{c}_1
\end{equation}
where the Jacobian of the monomials ($\partial \boldsymbol{c}_m / \partial \boldsymbol{c}_1$), thus necessary for computations on the tangent spaces of $\mathcal{C}^{(n, m)}$, is an analytic function of the monomials, and hence it is itself a linear function of $\boldsymbol{c}_m$. 

Eq. \eqref{eq :TSE def} is valid only in the vicinity of $\bar{\boldsymbol{x}}$, and converges as $m \rightarrow \infty$. The truncation error can be defined as:
\begin{equation}
    e_m = \lVert f(\boldsymbol{x}) - \Psi_m\boldsymbol{c}_m\rVert
\end{equation}
The value of this error increases with the deviation $\boldsymbol{c}_1$ from the reference.  In the case where the expanded function stands for the flow of a dynamic system, this error grows not only with larger deviations but also with increasing propagation time. To ensure that each temporal point of the trajectory is contained within an instantaneous domain of validity (a local neighborhood of the reference where the expansion is sufficiently valid), the larger deviation of $\boldsymbol{c}_1$ should be such that the truncation error at the final time of propagation is less than a chosen threshold $\varepsilon$. Thus, it is possible to define a region of validity $C_{\varepsilon}$ for $\boldsymbol{c}_1$ as follows:
\begin{equation}
    C_{\varepsilon} = \{\boldsymbol{c}_1: \lVert \boldsymbol{c}_1 \rVert \leq d_{\text{crit}}(\varepsilon)\}
\end{equation}
where $d_{\text{crit}}$ is defined as the largest deviation from the reference at $t_0$ such that the truncation error at $t_f$ is equal to the threshold $\varepsilon$. The value of $d_{\text{crit}}$ is determined \textit{a posteriori} once $\Psi_m$ has been computed. The process is a numerical iterative procedure that progressively increases the initial deviation along the direction of maximum perturbation, identified by the eigenvector associated with the largest eigenvalue of the Cauchy-Green tensor \cite{short_stretching_2015}, until the truncation error exceeds the chosen threshold.

\subsection{Differential Algebra}
Differential algebra is a mathematical approach that can be used to compute function derivatives on a computer, allowing easy computation of Taylor series expansions. Originally developed in the field of particle beam physics \cite{berz_chapter_1999}, DA has since found applications in space engineering, including continuous control \cite{di_lizia_high_2014, di_lizia_high_2014-1} and uncertainty propagation \cite{valli_nonlinear_2013, wittig_propagation_2015}. DA enables the computation of derivatives of any function $\boldsymbol{f}$ of $n$ variables up to arbitrary order $m$, alongside the function evaluation. This has an important consequence when the numerical integration of an ODE is performed using any integration scheme, since it relies on algebraic operations \cite{cellier_continuous_2006}. Using differential algebra, it is possible to compute an arbitrary order expansion of the flow of a general ODE with respect to the initial condition. \citet{gimeno_numerical_2023} demonstrate that applying methods such as Runge-Kutta, Taylor, or multistep to variational equations up to order $m$ yields the same results as performing the integration within the {differential algebra approach} using the same scheme.
The main advantage of this DA-based approach is that there is no need to integrate variational equations \cite{park_nonlinear_2006} to obtain high-order flow expansions. Therefore, by exploiting the DA framework, it is possible to build the matrix $\Psi_m$ at each integration node, which collects all the coefficients of the TSE,  easily and with limited computational effort. 
 In the Python environment, selected as the language in this work, two different free packages are available for DA: the ESA's PyAudi \cite{pyaudi-repository}
and DACEyPy \cite{daceypy-repository} 
(a Python front-end for the Differential Algebra Computational Engine DACE).

\section{Dynamical Model} \label{sec:dyn model}
Let us consider a rendezvous problem with a non-cooperative target. For the sake of simplicity, the dynamic models adopted in this work neglect effects such as solar radiation pressure and the gravitational attraction of the Sun, focusing solely on the primary gravitational influences established within the cislunar region. The target spacecraft is assumed to be on an NRHO belonging to the family of $L_2$ southern Halo orbits. The orbits in this family are of significant importance because they have been chosen as the Gateway's operational orbit \cite{lee_white_2019}. 
The dimensionless initial conditions used in this study to propagate the target's orbit are listed in Table \ref{tab: CR3BP IC}. Table \ref{tab: CR3BP IC} also includes the NRHO’s period ($T$) and Jacobi constant ($C$). The corresponding dimensionless units are reported in Table \ref{tab: CR3BP PAR}.

\begin{table}[h]
\centering
\begin{minipage}{0.45\textwidth}
\centering
        \caption{Dimensionless units.}
        \begin{tabular}{c c}
        \toprule
        \textbf{Parameter} & \textbf{Value} \\
        \midrule \midrule
        $x_0$ (-) & $1.01865930$ \\
        $y_0$ (-) & $0.0$ \\
        $z_0$ (-) & $-0.17967210$ \\
        $\dot{x}_0$ (-) & $8.74222438 \times 10^{-14}$ \\
        $\dot{y}_0$ (-) & $-0.09581408$ \\
        $\dot{z}_0$ (-) & $1.31415366 \times 10^{-12}$ \\
        $T$ (\si{days}) & $6.52499502$ \\
        $C$ (-) & $3.04997282$ \\
        \bottomrule
        \end{tabular}
        \label{tab: CR3BP IC}
\end{minipage}
\hfill
\begin{minipage}{0.45\textwidth}
\centering
        \caption{Earth–Moon NRHO parameters.}
            \begin{tabular}{c c}
        \toprule
        \textbf{Parameter} & \textbf{Value} \\
        \midrule \midrule
        Lenght unit LU (\unit{km}) & $389703$ \\
        Time unit TU (\unit{s}) & $382981$ \\
        Velocity unit VU (\unit{km/s}) & $1.017551785$ \\
        \bottomrule
        \end{tabular}
        \label{tab: CR3BP PAR}
\end{minipage}
\end{table}

Several applications of relative motion for NRHOs exist in the literature \cite{lizy-destrez_rendezvous_2019, scorsoglio_relative_2023, franzini_relative_2019}.  However, to implement the proposed guidance method with complex dynamics, this work adopts a fully nonlinear synodic relative equation of motion expressed in Cartesian coordinates. Both the target and the chaser are in the cislunar domain under the CR3BP assumption \cite{szebehely_chapter_1967}, therefore, the equations governing the dynamics of both spacecraft in the synodic frame are:
\begin{subequations}
\begin{align}
        \boldsymbol{\Ddot{r}}_\text{t} = -2\boldsymbol{\omega} \times \boldsymbol{\dot{r}}_\text{t} - \boldsymbol{\omega} \times (\boldsymbol{\omega} \times \boldsymbol{r}_\text{t}) - (1-\mu) \dfrac{\boldsymbol{r}_{1\text{t}}}{\lVert \boldsymbol{r}_{1\text{t}} \rVert^3} - \mu \dfrac{\boldsymbol{r}_{2\text{t}}}{\lVert \boldsymbol{r}_{2\text{t}} \rVert^3} \label{eq:rel eom 1}\\
            \boldsymbol{\Ddot{r}}_\text{c} = -2\boldsymbol{\omega} \times \boldsymbol{\dot{r}}_\text{c} - \boldsymbol{\omega} \times (\boldsymbol{\omega} \times \boldsymbol{r}_\text{c}) - (1-\mu) \dfrac{\boldsymbol{r}_{1\text{c}}}{\lVert \boldsymbol{r}_{1\text{c}} \rVert^3} - \mu \dfrac{\boldsymbol{r}_{2\text{c}}}{\lVert \boldsymbol{r}_{2\text{c}} \rVert^3} \label{eq:rel eom 2}
\end{align}
\end{subequations}
Here,  the subscript t indicates the target's quantities and c the chaser's ones. Let us define the synodic relative state vector as:
\begin{equation}
    \boldsymbol{x}  = \boldsymbol{x}_\text{c} - \boldsymbol{x}_\text{t} = [x , y , z , \dot{x} , \dot{y} , \dot{z} ]^\top = [\boldsymbol{\rho}^\top, \dot{\boldsymbol{\rho}}^\top]^\top 
\end{equation}
where $\boldsymbol{x}_\text{c} = [\boldsymbol{r}_\text{c}^\top, \dot{\boldsymbol{r}}_\text{c}^\top]^\top$ and $\boldsymbol{x}_\text{t} = [\boldsymbol{r}_\text{t}^\top, \dot{\boldsymbol{r}}_\text{t}^\top]^\top$  represent the chaser and target spacecraft state vectors in the synodic reference frame, respectively. Here, $\boldsymbol{\rho} = \boldsymbol{r}_\text{c} - \boldsymbol{r}_\text{t}$ denotes the relative position between the two spacecraft. The chaser's relative positions from the two primaries, $\boldsymbol{r}_{1\text{c}}$ and $\boldsymbol{r}_{2\text{c}}$, can be reformulated as:
   \begin{equation} \label{eq: r1 r2 rho rela}
    \begin{aligned}
        \boldsymbol{r}_{1\text{c}} = \boldsymbol{r}_{1\text{t}} + \boldsymbol{\rho} \\
        \boldsymbol{r}_{2\text{c}} = \boldsymbol{r}_{2\text{t}} + \boldsymbol{\rho} 
    \end{aligned}
    \end{equation}
By subtracting the adimensional Eq. \eqref{eq:rel eom 1} from Eq. \eqref{eq:rel eom 2}, and applying the relationship in Eq. \eqref{eq: r1 r2 rho rela}, the adimensional equations of motion in the relative synodic frame are obtained:
\begin{equation}
\label{eq:relative_cr3bp}
    \begin{cases}
    \Ddot{x}  = 2\dot{y}  + x  + (1-\mu)\left[\dfrac{\mu + x_\text{t}}{\lVert \boldsymbol{r}_{1\text{t}} \rVert ^3} - \dfrac{\mu + x_\text{t} + x }{\lVert \boldsymbol{r}_{1\text{t}} + \boldsymbol{\rho}\rVert ^3}\right] + \mu\left[\dfrac{\mu + x_\text{t}-1}{\lVert \boldsymbol{r}_{2\text{t}} \rVert ^3} - \dfrac{\mu + x_\text{t} + x -1}{\lVert \boldsymbol{r}_{2\text{t}} + \boldsymbol{\rho}\rVert ^3}\right]\\
    \Ddot{y}  = -2\dot{x}  + y  + (1-\mu)\left[\dfrac{y_\text{t}}{\lVert \boldsymbol{r}_{1\text{t}} \rVert ^3} - \dfrac{y_\text{t} + y }{\lVert \boldsymbol{r}_{1\text{t}} + \boldsymbol{\rho}\rVert ^3}\right] + \mu\left[\dfrac{y_\text{t}}{\lVert \boldsymbol{r}_{2\text{t}} \rVert ^3} - \dfrac{y_\text{t} + y }{\lVert \boldsymbol{r}_{2\text{t}} + \boldsymbol{\rho}\rVert ^3}\right]\\
    \Ddot{z}  = (1-\mu)\left[\dfrac{z_\text{t}}{\lVert \boldsymbol{r}_{1\text{t}} \rVert ^3} - \dfrac{z_\text{t} + z }{\lVert \boldsymbol{r}_{1\text{t}} + \boldsymbol{\rho}\rVert ^3}\right] + \mu\left[\dfrac{z_\text{t}}{\lVert \boldsymbol{r}_{2\text{t}} \rVert ^3} - \dfrac{z_\text{t} + z }{\lVert \boldsymbol{r}_{2\text{t}} + \boldsymbol{\rho}\rVert ^3}\right]\\
\end{cases} \hspace{-4em}
\end{equation}
The resulting dynamics can be compactly expressed in a state-space formulation as:
\begin{equation}
\label{eq:compact_ocp_dynamics}
    \dot{\boldsymbol{x}}  = \boldsymbol{f}(\boldsymbol{x} , t)
\end{equation}
where $\boldsymbol{f}: \mathbb{R}^7 \rightarrow \mathbb{R}^6$ represents the time-dependent dynamics function, which governs the evolution of the relative state vector $\boldsymbol{x} $.

\section{Optimal Control Problem} \label{sec:ocp}
\subsection{Statement of the problem}
An impulsive thrust model is considered for the chaser spacecraft, where impulsive maneuvers are performed at predefined points on a time grid.  The trajectory is discretized over $N+1$ points between the initial time $t_0$ and the final time $t_f$,
\begin{equation}
    t_0 < t_1 < \dots < t_{N-1} < t_N = t_f
\end{equation}
Each grid node allows for the possibility of an impulsive maneuver. Between two consecutive nodes, the spacecraft follows the natural relative dynamics as expressed in Eq. \eqref{eq:compact_ocp_dynamics}.  The relative initial condition of the chaser spacecraft is {$\boldsymbol{x}_{0} = [\boldsymbol{\rho}_{0}, \dot{\boldsymbol{\rho}}_{0}]^\top$}, and it must reach a specified relative final state  {$\boldsymbol{x}_{f} = [\boldsymbol{\rho}_{f}, \dot{\boldsymbol{\rho}}_{f}]^\top$}. The objective is to minimize the total propellant consumption, corresponding to minimizing the sum of the magnitudes of the velocity impulses at each node. The resulting OCP can be stated as:
\begin{subequations}
\label{eq:OCP}
\begin{align}
   &  \min_{\Delta\boldsymbol{v}_i} \; \sum_{i=1}^N \lVert \Delta\boldsymbol{v}_i\rVert  \; \; \text{s.t.} \\
   & 
   \begin{cases}
   \label{eq:cond_OCP}
       \boldsymbol{x}(t_0) = \boldsymbol{x}_0 \\
       \boldsymbol{x}(t_f) = \boldsymbol{x}_f \\
       \boldsymbol{x}(t_{i+1}^-) = \boldsymbol{\phi}(t_i, t_{i+1}, \boldsymbol{x}(t_i^+)
       ) \\
       \boldsymbol{r}(t_i^+) = \boldsymbol{r}(t_i^-) \\
       \boldsymbol{v}(t_i^+) = \boldsymbol{v}(t_i^-) + \Delta \boldsymbol{v}_i
   \end{cases} 
   \quad \; i=1, \dots, N
\end{align}
\end{subequations}
Here, time is discretized a priori and is not an optimization variable, although in the final part of this section, we propose a way to generalize the method to solve for free-final time rendezvous. Subscript $i$ is used to identify quantities at time $t_i$, whereas superscript $-$ and $+$ are used to distinguish the spacecraft state immediately before and after the impulse, respectively. The term $\boldsymbol{\phi}(t_i, t_{i+1}, \boldsymbol{x}(t_i))$ represents the flow of Eq. \eqref{eq:compact_ocp_dynamics} from $t_i$ to $t_{i+1}$ with initial condition $\boldsymbol{x}(t_i)$. The first and second conditions in Eq. \eqref{eq:cond_OCP} are the initial and final state boundary conditions, the third constraint imposes natural dynamics between maneuvers, the fourth guarantees position continuity through an instantaneous maneuver, while the last defines the impulsive velocity change.

\subsection{Formulation of the Optimization Problem in Monomial Coordinates}
\begin{figure}[h]
    \centering
    \includegraphics[width=0.9\linewidth]{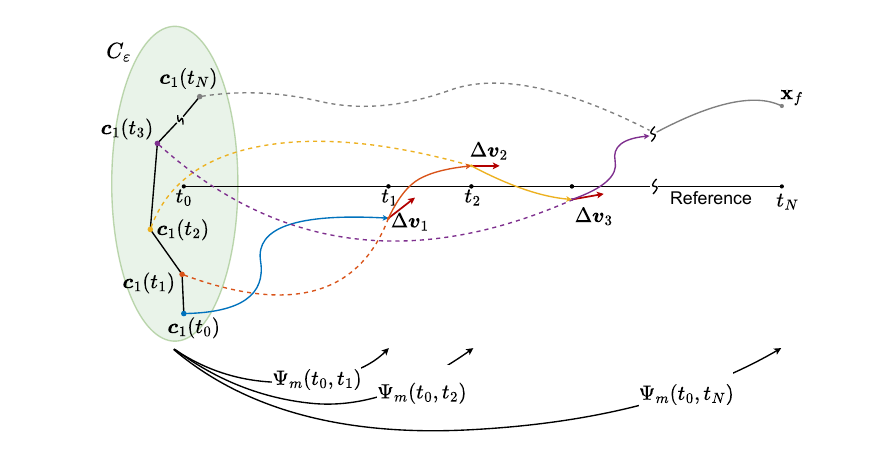}
    \caption{Impulsive maneuver trajectory transcription via $\Psi_m$.}
    \label{fig:impulsice_scheme_trasn}
\end{figure}

The optimization problem in Eq. \eqref{eq:OCP} is reformulated in terms of monomial coordinates. The spacecraft state at an epoch $t_i$ is approximated as:
\begin{equation}
    \boldsymbol{x}(t_i) \approx \Psi_m(t_0, t_i) \boldsymbol{c}_m(t_i)
\end{equation}
where $\Psi_m(t_0, t_i)$ is the $m$-th order monomial matrix of the flow of Eq. \eqref{eq:compact_ocp_dynamics}, which maps an initial condition at $t_0$ into the state at $t_i$, and $\boldsymbol{c}_m(t_i)$ is the vector of monomial coordinates of the state at $t_i$. The matrix $\Psi_m$ is computed using DA 
with the target orbit as reference.  As the target orbit is the reference of the expansion, the linear part of $\boldsymbol{c}_m(t_i)$ can be physically interpreted as the "fictitious" relative initial condition. 
By expressing the problem in monomial parametrization, the impulsive maneuver trajectory optimization is transformed into a path-planning problem \cite{burnett_rapid_2025}. At each node $t_i$, the chaser position, computed as $\Psi_m(t_0, t_i)\boldsymbol{c}_m(t_i)$, must match the position generated from  $\Psi_m(t_0, t_i)\boldsymbol{c}_m(t_{i-1})$, to ensure position continuity. However, this does not hold for velocity; the change in velocity determines the $i$-th impulse $\Delta\boldsymbol{v}_i$. This is illustrated in Figure \ref{fig:impulsice_scheme_trasn}, where the horizontal line represents the reference trajectory for the matrix $\Psi_m$. 
Starting from the initial condition $\boldsymbol{c}_1(t_0) = \boldsymbol{x}_0$, the chaser reaches $t_1$, where the first impulse is applied, shifting it to the trajectory generated by $\boldsymbol{c}_m(t_1)$. At the second node $t_2$, a second impulse is executed and the chaser is moved to the trajectory starting from $\boldsymbol{c}_m(t_2)$. This process continues until the chaser reaches the last node matching the target state $\boldsymbol{x}_f$. The figure shows the actual trajectory as a solid line, while the dashed line indicates the hypothetical trajectory that the chaser would have if it started from $\boldsymbol{c}_m(t_i)$ at $t_0$ under purely uncontrolled dynamics. To ensure accuracy in the approximation, each $\boldsymbol{c}_1(t_i)$ must lie in the region of validity $C_\varepsilon$, shown green in the figure. Consequently, because the entire relative state at any time $t_i$ is just governed by the value of $\boldsymbol{c}_m(t_i)$ or equivalently $\boldsymbol{c}_1(t_i)$, the trajectory optimization problem is reduced to the selection of the elements $\boldsymbol{c}_1(t_i)$ in $C_\varepsilon$.

\noindent At any given node $i$, the velocity impulse $\Delta\boldsymbol{v}_i$, can be reformulated as:
\begin{equation}
    \Delta\boldsymbol{v}_i = \boldsymbol{v}(t_i^+) - \boldsymbol{v}(t_{i}^-) = 
    \Psi_{v, m}(t_0, t_i) (\boldsymbol{c}_m(t_i) - \boldsymbol{c}_m(t_{i-1}))
\end{equation}
where $\Psi_{v, m}(t_0, t_i)$ represents the half of the monomial matrix corresponding to the velocity states. Similarly, the continuity condition for the position at node $i$ can be written as:
\begin{equation}
    \boldsymbol{r}(t_i^+) - \boldsymbol{r}(t_i^-) = 
    \Psi_{r,m}(t_0, t_i) (\boldsymbol{c}_m(t_i) - \boldsymbol{c}_m(t_{i-1})) = \boldsymbol{0}
\end{equation}
where now $\Psi_{r,m}(t_0, t_i)$ is the half of the monomial matrix for the position states. 
The optimization problem can now be reformulated in terms of the monomial coordinates, ensuring they remain within the validity region $C_\varepsilon$, as:
\begin{subequations}
\label{eq:ocp_stt}
\begin{align} \label{eq: nonLinearCost}
   &  \min_{\boldsymbol{c}_1(t_i)} \; \sum_{i=1}^N \lVert \Psi_{v,m}(t_0, t_i) (\boldsymbol{c}_m(t_i) - \boldsymbol{c}_m(t_{i-1}))\rVert \\
   & \text{s.t.} \quad
   \begin{cases}
   \label{eq:constraint_ocp}
       \boldsymbol{c}_m(t_0) = \boldsymbol{c}_{m, 0}=\boldsymbol{E}_m(\boldsymbol{x}_0) \\
       \Psi_{m}(t_0, t_N)\boldsymbol{c}_m(t_N) = \boldsymbol{x}_f \\
       \boldsymbol{c}_m(t_i) = \boldsymbol{E}_m(\boldsymbol{c}_1(t_i))\\
       \Psi_{r,m}(t_0, t_i) (\boldsymbol{c}_m(t_i) - \boldsymbol{c}_m(t_{i-1})) = \boldsymbol{0} \\
       \boldsymbol{c}_1(t_i) \in C_{\varepsilon}
   \end{cases} 
   \quad  \; i=1, \dots, N
\end{align}
\end{subequations}
The third constraint in Eq. \eqref{eq:constraint_ocp}, which imposes the relationship between the linear and nonlinear parts of the set of monomials, introduces non-convexity into the problem. To address this, a sequential convex programming approach is proposed, which iteratively solves a convex approximation of the problem. 

\subsection{Sequential Convex Programming}
Following a standard workflow of the sequential convex programming method \cite{malyuta_convex_2022}, the non-convex terms are linearized. The linearization of $\boldsymbol{c}_m(t_i)$ results in:
\begin{equation}
    \boldsymbol{c}_m(t_i) \approx \boldsymbol{c}_m^{'}(t_i) + \dfrac{\partial \boldsymbol{c}_m}{\partial \boldsymbol{c}_1}\bigg|_{i} \delta\boldsymbol{c}_1(t_i)
\end{equation}
where the superscript $\square^{'}$ denotes the reference value. 
According to the third constraint in Eq. \eqref{eq:constraint_ocp}, every monomial $\boldsymbol{c}_m(t_i)$ must lie on the manifold $\mathcal{C}^{(n, m)}$. Therefore, the underlying mathematical idea is to express each monomial coordinate $\boldsymbol{c}_m(t_i)$ as a reference value $\boldsymbol{c}_m^{'}(t_i)$ on the manifold plus a deflection in the plane that is locally tangent to the manifold at that reference value.

Substituting the linearization of the $\boldsymbol{c}_m(t_i)$ into Eq. \eqref{eq:ocp_stt}, the resulting convex sub-problem is expressed as:
\begin{subequations}
\label{eq:scp_ocp}
\begin{align}
   \label{eq:scp_ocp_obj}
   \min_{\delta\boldsymbol{c}_1(t_i), \boldsymbol{s}_i, \boldsymbol{s}_{end}} \; 
   \begin{array}{cc}
       \left \lVert \Psi_{v, m}(t_0, t_1)\left(\boldsymbol{c}_m^{'}(t_1) + \dfrac{\partial \boldsymbol{c}_m}{\partial \boldsymbol{c}_1}\bigg|_{1} \delta\boldsymbol{c}_1(t_1) - \boldsymbol{c}_m(t_0) \right) \right\rVert + w\left(\sum_{i=1}^N\lVert\boldsymbol{s}_i \rVert^2 + \lVert\boldsymbol{s}_{end}\rVert^2 \right) \\
       +\sum_{i=2}^N\left \lVert \Psi_{v, m}(t_0, t_i)\left(\boldsymbol{c}_m^{'}(t_i) + \dfrac{\partial \boldsymbol{c}_m}{\partial \boldsymbol{c}_1}\bigg|_{i} \delta\boldsymbol{c}_1(t_i) - \boldsymbol{c}_m^{'}(t_{i-1}) - \dfrac{\partial \boldsymbol{c}_m}{\partial \boldsymbol{c}_1}\bigg|_{i-1} \delta\boldsymbol{c}_1(t_{i-1}) \right) \right\rVert  
   \end{array}
   \end{align}
   \vspace{-0.7cm}
   \begin{align}
   \text{s.t.} \quad
   \label{eq:scp_ocp_cons}
   \begin{cases} 
       \Psi_{m}(t_0, t_N)\left(\boldsymbol{c}_m^{'}(t_N) + \dfrac{\partial \boldsymbol{c}_m}{\partial \boldsymbol{c}_1}\bigg|_{N} \delta\boldsymbol{c}_1(t_N)\right) + \boldsymbol{s}_{end} = \boldsymbol{x}_f \\
        \Psi_{r, m}(t_0, t_1)\left(\boldsymbol{c}_m^{'}(t_1) + \dfrac{\partial \boldsymbol{c}_m}{\partial \boldsymbol{c}_1}\bigg|_{1} \delta\boldsymbol{c}_1(t_1) - \boldsymbol{c}_m(t_0) \right) + \boldsymbol{s}_1= \boldsymbol{0} \\
        \Psi_{r, m}(t_0, t_i)\left(\boldsymbol{c}_m^{'}(t_i) + \dfrac{\partial \boldsymbol{c}_m}{\partial \boldsymbol{c}_1}\bigg|_{i} \delta\boldsymbol{c}_1(t_i) - \boldsymbol{c}_m^{'}(t_{i-1}) - \dfrac{\partial \boldsymbol{c}_m}{\partial \boldsymbol{c}_1}\bigg|_{i-1} \delta\boldsymbol{c}_1(t_{i-1}) \right) + \boldsymbol{s}_i = \boldsymbol{0} \\
        \hspace{10.37cm}  \; i=2, \dots, N \\
       \lVert \boldsymbol{c}_1^{'}(t_i) + \delta\boldsymbol{c}_1(t_i)\rVert \leq {d_{\text{crit}}(\varepsilon)} \quad \quad  \; i=1, \dots, N 
    \end{cases} \hspace{-3.5em}
\end{align}
\end{subequations}
This convex sub-problem involves not only the substitution of the local tangent plane approximation but also the introduction of new variables. The variables $\boldsymbol{s}_i$ with $i=0, \dots, N $ and $\boldsymbol{s}_{end}$ are the slack variables introduced to prevent artificial infeasibility that may arise due to the linearization.  All slack variables are penalized in the problem's objective by a positive scalar weight $w$. The free variable of the convex subproblem consists of the tangent-plane approximation variations $\delta\boldsymbol{c}_1(t_i)$ and the slack variables. To avoid the issue of artificial unboundedness, a trust region constraint is enforced only on the tangent-plane variations:
\begin{equation}
    \lVert (\delta\boldsymbol{c}_1(t_1)^\top, \dots, \delta\boldsymbol{c}_1(t_N)^\top)^\top \rVert \leq \eta^{(k)}
\end{equation}
The trust region radius ($\eta$) is updated between iterations with the following rule \cite{malyuta_convex_2022}:
\begin{equation} \label{eq: eta update}
    \eta^{(k+1)} = 
    \begin{cases}
        \beta \eta^{(k)} & \text{if} \; \zeta_0 \leq \zeta^{(k)} < \zeta_1\\      
       \eta^{(k)} &  \text{if} \; \zeta_1 \leq \zeta^{(k)} \zeta_2 \\      
        \alpha \eta^{(k)}&\text{if} \;  \zeta^{(k)} > \zeta_2\\ 
    \end{cases}
\end{equation}
where  $\beta < 1$ and $\alpha > 1$ are the trust region shrinking and growing rate, respectively, and $\zeta$ is the index of convexification accuracy at iteration $k$ of the SCP defined as \cite{malyuta_convex_2022}:
\begin{equation}
    \zeta = \dfrac{\text{actual cost improvement}}{\text{predicted cost improvement}} 
\end{equation}
Here, the actual cost improvement is the difference between the nonlinear cost evaluated at the reference monomials $\boldsymbol{c}_m^{'}(t_i)$  and the nonlinear cost at the new subproblem solution $\boldsymbol{c}_m(t_i)$. The predicted cost improvement is the difference between the nonlinear cost at the reference and the linearized cost. The thresholds $\zeta_0, \zeta_1, \zeta_2 \in (0, 1)$ are the smallest, intermediate, and largest acceptance metrics, respectively.  

After each SCP iteration, the solution from the convex sub-problem needs to be projected back onto the manifold. The tangent-plane approximations of the variations $\delta\boldsymbol{c}_1(t_i)$ are used to update the linear part of $\boldsymbol{c}_m^{'}(t_i)$ and by Eq. \eqref{eq:mon_pro} the new vector of monomial coordinates is computed.
The flowchart of the nonlinear convex method for impulsive guidance is outlined in Figure \ref{fig:scp_flowchart}. The number of impulses is not defined a priori but emerges as part of the converged guidance solution. Once the SCP algorithm converges, the resulting monomials $\boldsymbol{c}_m(t_i)$ undergo a post-processing phase, during which the sequence of impulsive maneuvers is identified. The linear part of the monomial at node $i$ is compared with the linear part of the previous node.  If the difference between them exceeds a specific threshold ($\delta_c$),  an impulsive maneuver is scheduled to occur at time $t_i$ with $\Delta \boldsymbol{v}_i = \Psi_{v, m}(t_0, t_i) (\boldsymbol{c}_m(t_i) - \boldsymbol{c}_m(t_{i-1}))$; otherwise, no maneuver is executed. In the rendezvous examples in work, $\delta_c$ is fixed at $10^{-4}$. 

\begin{figure}[h]
\centering

    \tikzset{
  block/.style={
    draw,
    rectangle,
    rounded corners,
    align=center,
    text width=5cm, 
    minimum height=1.2cm
  },
  decision/.style={
    draw,
    diamond,
    aspect=2,
    align=center,
    text width=3cm, 
    inner sep=1pt
  },
  arrow/.style={-{Stealth}, thick}
}
\resizebox{\textwidth}{!}{%
\begin{tikzpicture}[node distance=1cm]

\node[block] (init) {
 Set initial state $\boldsymbol{x}_0$, target state $\boldsymbol{x}_f$, and SCP parameters (Table \ref{tab: SCP parameter})
};

\node[block, below=of init] (guess) {
  Generate an initial guess by solving the linear version of the guidance problem \eqref{eq:ocp_stt} setting $m=1$.
};

\node[block, right= of guess, text width=3cm] (iter) {
  Initialize iteration \\ counter $k = 0$
};

\node[decision, right= of iter] (cond1) {
  $k <  N_\text{max}$
};

\node[block, text width=7cm] (derivatives) at ([xshift=5cm,yshift=2cm]cond1.north east) {
  Compute $\dfrac{\partial \boldsymbol{c}_m}{\partial \boldsymbol{c}_1}\bigg|_{i}$ for the current guess $\boldsymbol{c}_1(t_i)$
};

\node[block, below=of derivatives, text width=7cm] (solve) {
   Solve the convex sub-problem (Eq.\eqref{eq:scp_ocp}) to find $\delta\boldsymbol{c}_1(t_i)$
};

\node[block, below=of solve, text width=7cm] (update) {
  $\boldsymbol{c}_1(t_i) \leftarrow
  \boldsymbol{c}_1(t_i) + \delta\boldsymbol{c}_1(t_i)$
};

\node[decision, below=1.5cm of update, text width=3cm] (optimal) {
  Optimal and feasible
};

\node[block, below=of cond1, left=1.4cm of optimal] (thrust) {
  Update thrust region using\\
  Eq.\eqref{eq: eta update}
};

\node[block, below=of cond1, text width=3cm] (update-trust) {
  Increment iteration counter: $k = k + 1$
};

\node[block, below=of optimal] (final) {
  Optimal solution
};

\node[block, left=of final, text width=7cm] (post) {
Post-process the SCP optimal solution to identify the impulsive maneuvers.
};

\draw[arrow] (init) -- (guess);
\draw[arrow] (guess) -- (iter);
\draw[arrow] (iter) -- (cond1);

\draw[arrow] (cond1.north) |- node[above]{yes} (derivatives.west);
\draw[arrow] (derivatives.south) -- (solve.north);
\draw[arrow] (solve) -- (update);
\draw[arrow] (update) -- (optimal);

\draw[arrow] (optimal) -- node[right]{yes} (final);
\draw[arrow] (optimal.west) -- node[above]{no} (thrust.east);
\draw[arrow] (thrust.north) -- (update-trust.south);
\draw[arrow] (update-trust.north) -- (cond1.south);

\draw[arrow] (final.west) -- (post.east);

\end{tikzpicture}%
}
\caption{Flowchart of the fixed-final time nonlinear convex impulsive guidance.}
\label{fig:scp_flowchart}
\end{figure}

\subsection{Monomial Method for Closed-Loop On-board Guidance}
\begin{figure}[h]
    \centering
    \includegraphics[width=0.7\linewidth]{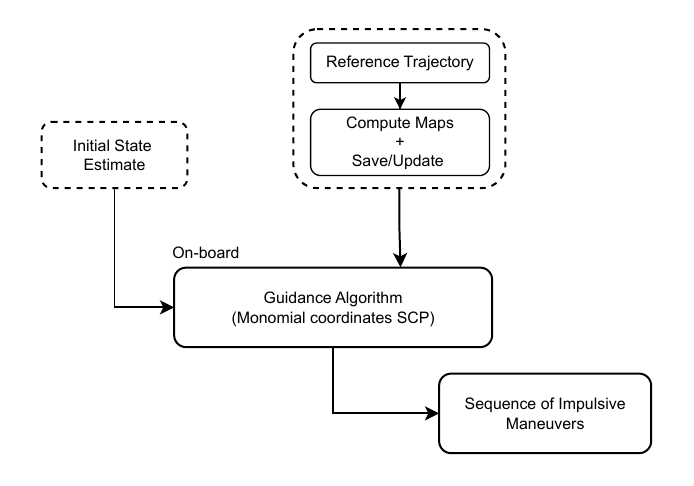}
    \caption{Block diagram illustrating the closed-loop guidance.}
    \label{fig:bd closed-loop}
\end{figure}

The monomial-based guidance method provides a robust and computationally efficient guidance method suitable for closed-loop on-board implementation on resource-constrained platforms. At each guidance cycle, the algorithm receives the estimated relative state from the navigation system and accesses the precomputed Taylor expansion map $\boldsymbol{\Psi}_m$ of the dynamics, which is stored in onboard memory. In this framework, deviations from the nominal trajectory arising from unmodeled perturbations or execution uncertainties are corrected through successive re-optimizations over a receding horizon time of flight. Figure \ref{fig:bd closed-loop} illustrates the block-diagram of the closed-loop guidance scheme. 

The most computationally expensive operation in the proposed framework is the computation of the dynamical map $\boldsymbol{\Psi}_m$. However, the map is not regenerated every time the guidance solution is executed. Instead, once the map has been computed for a given target orbit, it can be reused across multiple guidance cycles. Assume that the Taylor expansion map $\boldsymbol{\Psi}_m$ has been computed over the interval  $[\tau_0, \tau_f]$, while a navigation update is received at a generic time  $\tau_k \in [\tau_0, \tau_f]$. Because the guidance formulation with monomial method requires a dynamical map originating from the current time $\tau_k$, a direct recomputation of the map would appear necessary. However, by exploiting the composition properties of the dynamical map, the required map can be obtained through algebraic manipulation of the precomputed one. In \cite{park_nonlinear_2007}, the rules for STT composition are derived. For instance, in the case of a first-order map, the following relation holds:
\begin{equation}
    \boldsymbol{\Psi}_1(\tau_k, \tau_s) = \boldsymbol{\Psi}_1(\tau_0, \tau_s) \boldsymbol{\Psi}_1(\tau_0, \tau_k)^{-1}
\end{equation}
with $\tau_r \in [\tau_0, \tau_f]$. This property implies that, once the map has been generated along the reference trajectory, the map from an arbitrary point in time to some future instant can be obtained through simple algebraic operations. Therefore, the guidance algorithm can promptly incorporate navigation updates without additional integration or map recomputation.
Nevertheless, the map may become invalid if the relative distance between the chaser and the target grows beyond the region of validity of the Taylor expansion or if the target orbit is modified. In these situations, the map must be recomputed, either on the ground and uploaded to the spacecraft or autonomously onboard between maneuver executions.

\subsection{Accommodating State Constraints} \label{sec:scp cons}
In the guidance problem, let us now introduce two sets of time-dependent constraints: minimum-range and approach cone. 
Previous work has shown that the minimum-range constraint can be reformulated using a Taylor map \cite{burnett_rapid_2025}.
Motivated by this view, in this work the monomial parameterization is used to formulate another constraint: the approach cone. These two constraints are not active simultaneously: the minimum-range constraint is active for the first part of the flight, while the approach cone constraint is active in the final portion when the spacecraft inter-range separation is small. 

The instantaneous chaser-target range is defined as the Euclidean norm of their relative position: 
\begin{equation}
    \rho(t) = \lVert \boldsymbol{\rho}(t) \rVert = \sqrt{x(t)^2+y(t)^2+z(t)^2}
\end{equation}
The minimum-range constraint can be equivalently reformulated as $\rho(t)^2\geq \rho_\text{min}(t)^2$. This expression can be approximated using a TSE with respect to the origin. Define a scalar multivariable function $F:\mathbb{R}^7 \rightarrow \mathbb{R}$ such that 
\begin{equation}
    F(\boldsymbol{x}_0, t) = \lVert\boldsymbol{\phi}_r(t_0, \boldsymbol{x}_0; t) \rVert^2
\end{equation} 
where $\boldsymbol{x}_0$ is the initial condition of the relative motion at a time $t_0$ and $\boldsymbol{\phi}_r$ is the position-related flow of relative motion. This function is nonlinear in terms of relative position, but it is possible to exploit differential algebra to write the constraints as a linear function of $\boldsymbol{c}_m(t_i)$ as:
\begin{equation}
    F(\boldsymbol{x}_0, t_i)= \Gamma_m(t_0, t_i)\boldsymbol{c}_m(t_i)
\end{equation}
here $\Gamma_m(t_0, t_i)$ is the counterpart of $\Psi_m(t_0, t_i)$ for the function $ F(\boldsymbol{x}_0, t)$, and $\boldsymbol{c}_m$ is the monomial combination for $\boldsymbol{x}_0$. Therefore, the constraint can be reformulated as:
\begin{equation}
 \Gamma_m(t_0, t_i)\boldsymbol{c}_m(t_i) \geq \rho_\text{min}(t_i)^2 \quad  \; i=1, \dots, N_\text{wp}-1
\end{equation}
where  $N_\text{wp} < N$ defines the time $t_{N_\text{wp}}$ until which the minimum-range constraint is active. The expansion is performed around the reference trajectory, making the constant part of the expansion zero.

Next, consider the approach cone constraint. The unit vector $\boldsymbol{n}$ represents the docking axis, which is fixed in the target frame, while $\gamma$ is the semi-aperture angle of the approach cone, which typically ranges from {$10$ deg to $15$ deg}.  
The approach cone constraint can be imposed as \cite{lu_autonomous_2013}:
\begin{equation}
\label{eq:AC}
    \lVert \boldsymbol{r}(t) - \boldsymbol{r}_f \rVert \cos{\gamma} \leq \boldsymbol{n}^\top(\boldsymbol{r}(t) - \boldsymbol{r}_f)
\end{equation}

This formulation is convex, and according to the DCP rules \cite{diamond_cvxpy_2016}, its practical implementation into the SCP is straightforward. However, as with the minimum range, this constraint is expressed as a Taylor series {for convenience}. Since the left-hand side of the Eq. \eqref{eq:AC} is strictly positive, squaring both sides yields the identical formulation:  $\lVert \boldsymbol{r}(t) - \boldsymbol{r}_f \rVert^2 \cos^2{\gamma} \leq (\boldsymbol{n}^\top(\boldsymbol{r}(t) - \boldsymbol{r}_f))^2$. 
Let a scalar function $G:\mathbb{R}^7 \rightarrow \mathbb{R}$ be defined as
\begin{equation}
    G(\boldsymbol{x}_0, t) = \lVert \boldsymbol{\phi}_r(t_0, \boldsymbol{x}_0; t) - \boldsymbol{r}_f \rVert^2 \cos^2{\gamma} - (\boldsymbol{n}^\top(\boldsymbol{\phi}_r(t_0, \boldsymbol{x}_0; t) - \boldsymbol{r}_f))^2
\end{equation}
and via differential algebra we can impose $G$ as a linear function of $\boldsymbol{c}_m(t_i)$:
\begin{equation}
    G(\boldsymbol{x}_0, t_i) = G(\boldsymbol{0}, t_i) + \Theta_m(t_0, t_i)\boldsymbol{c}_m(t_i)
\end{equation}
This reformulation allows the constraint to be expressed in monomial coordinates expanded with respect to the origin, as follows:
\begin{equation}
    G(\boldsymbol{0}, t_i) +\Theta_m(t_0, t_i)\boldsymbol{c}_m(t_i) \leq 0 \quad  \; i=N_\text{wp}, \dots, N
\end{equation}

In order to apply these constraints to the SCP subproblem in Eq. \eqref{eq:scp_ocp},  they must be linearized following the same procedure exploited for the objective and the constraints. This leads to two {different form of} linear constraints, with additional slack variables introduced and penalized in the cost function as in Eq. \eqref{eq:scp_ocp_obj}. The linearized constraints are formulated as:
\begin{subequations}
    \begin{align}
        \Gamma_m(t_0, t_i)\left(\boldsymbol{c}_m^{'}(t_i) + \dfrac{\partial \boldsymbol{c}_m}{\partial \boldsymbol{c}_1}\bigg|_{i} \delta\boldsymbol{c}_1(t_i)\right) + s_{\text{range}, i} \geq \rho_{\text{min}}(t_i)^2 \quad \; i=1, \dots, N_\text{wp}-1 \label{eq:scp range} \\
        G(\boldsymbol{0}, t_i) + \Theta_m(t_0, t_i)\left(\boldsymbol{c}_m^{'}(t_i) + \dfrac{\partial \boldsymbol{c}_m}{\partial \boldsymbol{c}_1}\bigg|_{i} \delta\boldsymbol{c}_1(t_i)\right) + s_{\text{a.c.}, i} \leq 0 \quad \; i=N_\text{wp}, \dots, N \label{eq:scp cone}
    \end{align}
\end{subequations}
The presence of these two constraints adds $N$ new slack variables, making the optimization variable of higher dimension.

\subsection{Extension to Free-Time Rendezvous} \label{sec:free time theory}
Let us now explore how to generalize the proposed SCP-based method to handle a free-final time rendezvous, allowing the maneuver times $t_i$ to enter the optimization.  {Using} the formulation of the problem in Eq. \eqref{eq:ocp_stt} in terms of the monomial coordinates, we state the free-time problem as follows:
\begin{subequations}
\label{eq:ocp_stt_free}
\begin{align}
   &  \min_{\boldsymbol{c}_1(t_i),\, t_i} \; \sum_{i=1}^{N_l} \lVert \Psi_{v,m}(t_0, t_i) (\boldsymbol{c}_m(t_i) - \boldsymbol{c}_m(t_{i-1}))\rVert \\
   & \text{s.t.} \quad
   \begin{cases}
   \label{eq:constraint_ocp_free}
       \boldsymbol{c}_m(t_0) = \boldsymbol{c}_{m, 0} \\
       \Psi_{m}(t_0, t_N)\boldsymbol{c}_m(t_N) = \boldsymbol{x}_f \\
       \boldsymbol{c}_m(t_i) = \boldsymbol{E}_m(\boldsymbol{c}_1(t_i))\\
       \Psi_{r,m}(t_0, t_i) (\boldsymbol{c}_m(t_i) - \boldsymbol{c}_m(t_{i-1})) = \boldsymbol{0} \\
       \boldsymbol{c}_1(t_i) \in C_{\varepsilon} \\
       t_i \in [t_0, t_\text{max}] \\
       t_i > t_{i-1}
   \end{cases} 
   \quad \; i=1, \dots, N_l
\end{align}
\end{subequations}
The optimization variable now includes both the monomial coordinates $\boldsymbol{c}_1(t_i)$ and maneuver/burn times $t_i$ $\forall i \in \{1, N_l\}$ and two additional constraints are included in the problem which impose that each burn time lies within an "outer-bound" time range $[t_0, t_\text{max}]$ and that maneuver times are strictly increasing. Inspired by the fixed-time rendezvous results, where the optimal solution requires only a limited number of impulsive maneuvers, we let $N_l \ll N$ to ensure a reasonable optimization problem size. The outer-bound time interval is divided over $N+1$ nodes, with the monomial matrix $\Psi_m(t_0, t_i)$ available at each node. By interpolating $\Psi_m(t_0, t_i)$ across the nodes, e.g., via spline interpolation, we obtain a smooth approximation of $\Psi_m(t_0, t)$ over $[t_0, t_\text{max}]$. 

The third constraint in the problem is still non-convex, hence, as in the fixed-time guidance, we apply SCP techniques. The terms $\boldsymbol{c}_m(t_i)$ are linearized using tangent plane deviations, while the map $\Psi_m(t_0, t_i)$ is linearized with respect to a small change ($\delta t_i$) in the $i$-th burn time as:
 \begin{equation}
     \Psi_m(t_0, t_i^{'} + \delta t_i) = \Psi_m(t_0, t_i^{'}) + \dfrac{\partial \Psi_m(t_0, t)}{\partial t}\bigg|_{i} \delta t_i
 \end{equation}
where the superscript $\square^{'}$ denotes the results of the previous iteration, and the derivative terms are computed from the interpolating function. Substituting the linearized $\Psi_m$ into the convex sub-problem in Eq. \eqref{eq:scp_ocp}, and neglecting terms of order higher than one, yields the convex sub-problem for the free-time case. Like the fixed-time problem case, the convex sub-problem includes slack variables to avoid artificial infeasibility, and the free variables include now tangent plane deviations, infinitesimal time variations, and slacks. To avoid artificial unboundedness, we enforce trust region constraints:
\begin{equation}
\begin{cases}
        \lVert (\delta\boldsymbol{c}_1(t_1)^\top, \dots, \delta\boldsymbol{c}_1(t_{N_l})^\top)^\top \rVert \leq \eta^{(k)} \\
        \lVert ( \delta t_1, \dots, \delta t_{N_l})^\top \rVert \leq \eta^{(k)} \\
\end{cases}
\end{equation}
After each solution of the convex sub-problem, the infinitesimal time changes are used to update the maneuvering times while tangent-plane deviations define new $\boldsymbol{c}_1(t_i)$. The converged SCP solution inherits interpolation error from $\Psi_m$, in addition to the approximation error due to truncating at finite order $m$. Therefore, the solution of the SCP is post-processed by "snapping" each optimal burn time to the nearest node from the original grid of $N+1$ nodes. At this point, either $\Psi_m$-based differential correction approach (Stage $2$ in \cite{burnett_rapid_2025}) or the fixed-time nonlinear convex guidance scheme refines the solution, allowing maneuvers only at these selected nodes. Although this yields a sub-optimal solution, it removes the interpolation errors, leaving only the truncation ones. The SCP free-time flowchart is reported in Figure \ref{fig:scp_free_flowchart}.

\begin{figure}[h]
\centering

    \tikzset{
  block/.style={
    draw,
    rectangle,
    rounded corners,
    align=center,
    text width=5cm, 
    minimum height=1.2cm
  },
  decision/.style={
    draw,
    diamond,
    aspect=2,
    align=center,
    text width=3cm, 
    inner sep=1pt
  },
  arrow/.style={-{Stealth}, thick}
}
\resizebox{\textwidth}{!}{%
\begin{tikzpicture}[node distance=1cm]

\node[block] (init) {
 Set initial state $\boldsymbol{x}_0$, target state $\boldsymbol{x}_f$, and SCP parameters (Table \ref{tab: SCP parameter})
};

\node[block, below=of init] (time-init) {
  Set final time guess $t_f \leq t_\text{max}$.
};

\node[block, below=of time-init] (guess) {
  Generate an initial guess by solving the linear version of the guidance problem \eqref{eq:ocp_stt}, with $m=1$ and the initial guess for the final time.
};

\node[block, right= of guess, text width=3cm] (iter) {
  Initialize iteration \\ counter $k = 0$
};

\node[decision, right= of iter] (cond1) {
  $k <  N_\text{max}$
};

\node[block, text width=7cm] (derivatives) at ([xshift=5cm,yshift=2cm]cond1.north east) {
 Compute $\dfrac{\partial \boldsymbol{c}_m}{\partial \boldsymbol{c}_1}\bigg|_{i}$ and $\dfrac{\partial \Psi_m(t_0, t)}{\partial t}\bigg|_{i}$ for the current guess $\boldsymbol{c}_1(t_i)$ and $t_i$.
};

\node[block, below=of derivatives, text width=7cm] (solve) {
   Solve the free-time convex sub-problem to find $\delta\boldsymbol{c}_1(t_i)$ and $\delta t_i$
};

\node[block, below=of solve, text width=7cm] (update) {
  $\boldsymbol{c}_1(t_i) \leftarrow
  \boldsymbol{c}_1(t_i) + \delta\boldsymbol{c}_1(t_i)$
};

\node[block, below=of update, text width=7cm] (time-update) {
  $t_i \leftarrow t_i + \delta t_i$
};

\node[decision, below=1.5cm of time-update, text width=3cm] (optimal) {
  Optimal and feasible
};

\node[block, below=of cond1, left=1.4cm of optimal] (thrust) {
  Update thrust region using\\
  Eq.\eqref{eq: eta update}
};

\node[block, below=of cond1, text width=3cm] (update-trust) {
  Increment iteration counter: $k = k + 1$
};

\node[block, below=of optimal] (final) {
  Output optimal solution
};

\node[block, left=of final, text width=7cm] (post) {
Post-process the free-final time SCP optimal solution to remove interpolation error.
};

\node[block, left=of post, text width=7cm] (post-post) {
Post-process the corrected sub-optimal solution to identify the impulsive maneuvers.
};

\draw[arrow] (init) -- (time-init);
\draw[arrow] (time-init) -- (guess);
\draw[arrow] (guess) -- (iter);
\draw[arrow] (iter) -- (cond1);

\draw[arrow] (cond1.north) |- node[above]{yes} (derivatives.west);
\draw[arrow] (derivatives.south) -- (solve.north);
\draw[arrow] (solve) -- (update);
\draw[arrow] (update) -- (time-update);
\draw[arrow] (time-update) -- (optimal);

\draw[arrow] (optimal) -- node[right]{yes} (final);
\draw[arrow] (optimal.west) -- node[above]{no} (thrust.east);
\draw[arrow] (thrust.north) -- (update-trust.south);
\draw[arrow] (update-trust.north) -- (cond1.south);

\draw[arrow] (final.west) -- (post.east);
\draw[arrow] (post.west) -- (post-post.east);

\end{tikzpicture}%
}
\caption{Flowchart of the free-final time nonlinear convex algorithm.}
\label{fig:scp_free_flowchart}
\end{figure}

\section{Numerical Results} \label{sec:results}
This section presents the results obtained using the proposed nonlinear convex algorithm. All simulations for the test cases have been carried out on a Dell XPS $15$, equipped with an $11$-th Gen Intel(R) Core(TM) i$7$-$11800$H at $2.30$ GHz and with $16$ GB RAM. 
DACEyPy is used for automatic differentiation to compute the partials later used to build $\Psi_m$. In these examples, integrations of the dynamics are carried out using a Runge-Kutta 7/8 method with adaptive step size with relative and absolute tolerances set to $10^{-13}$ and $10^{-14}$, respectively. The convex subproblems are solved using the CLARABEL solver \cite{goulart_clarabel_2024}, with both relative and absolute tolerances set to $10^{-10}$.
\begin{table}[h]
    \centering 
    \caption{SCP parameters.}
    \begin{tabular}{ c l c }
    \toprule
     \textbf{Symbol} & \textbf{Quantity} & \textbf{Value}  \\
    \midrule \midrule
    $N_{\text{max}}$ &Maximum number of iteration& $25$ \\
    $\eta^{(0)}$& Initial trust radius & $0.005$    \\
    $\eta_{\text{max}}$ & Maximum trust radius & $0.5$   \\
    $\eta_{\text{min}}$ & Minimum trust radius & $5\e{-7}$   \\
    $\text{tol}_x$ & Minimum tolerance between the states of consecutive solutions & $5\e{-7}$   \\
    {$\zeta_0$} & Smallest acceptance metric & $0$  \\
    {$\zeta_1$} & Intermediate acceptance metric & $0.25$  \\
    {$\zeta_2$} & Largest acceptance metric & $0.90$   \\
    $\alpha$ & Trust region expansion factor  &    $2$  \\
    $\beta$ & Trust region reduction factor  &    $0.5$  \\
    \bottomrule
    \end{tabular}
    \label{tab: SCP parameter}
\end{table}
The parameters used in the SCP algorithm are summarized in Table \ref{tab: SCP parameter}. These values are constant across all simulations, except for the weight of the slack variable, $w$, which is tuned to ensure convergence and minimize guidance and open-loop errors.  With open-loop error ($\boldsymbol{e}_{\text{op}}$), we refer to the difference between the chaser's desired relative position and the actual position of the spacecraft after applying the impulsive maneuvers and it is present because the flow approximation has a truncation error due to the finite order of expansion $m$. While with $\boldsymbol{e}_\text{g}$ we indicate the guidance error, which is the difference between the desired relative position of the chaser and the position at the final time approximated by $\Psi_m$.

\begin{figure}[h]
    \centering
    \begin{minipage}{0.45\textwidth}
    \centering
    \includegraphics[width=\linewidth]{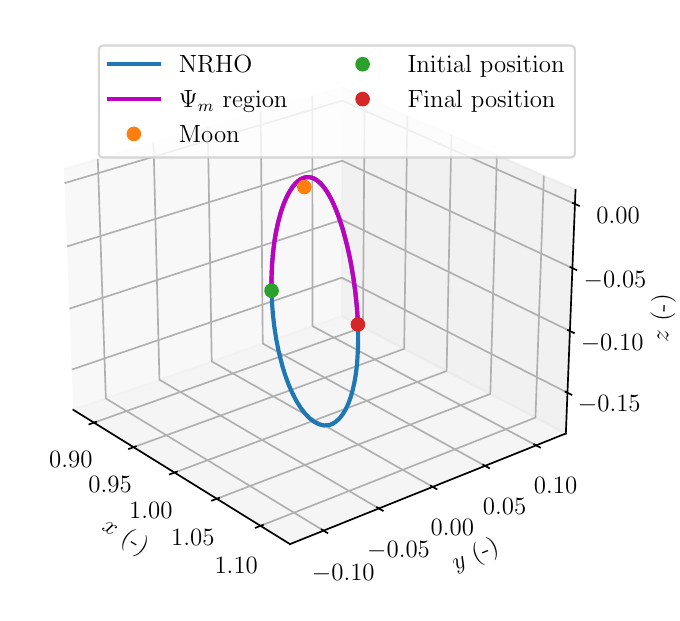}
    \caption{Target's NRHO and time interval of the generated $\Psi_m$ for  rendezvous {scenario} 1.}
    \label{fig:caseB_nrho}
    \end{minipage}
    \hfill
    \begin{minipage}{0.45\textwidth}
    \centering
    \includegraphics[width=\linewidth]{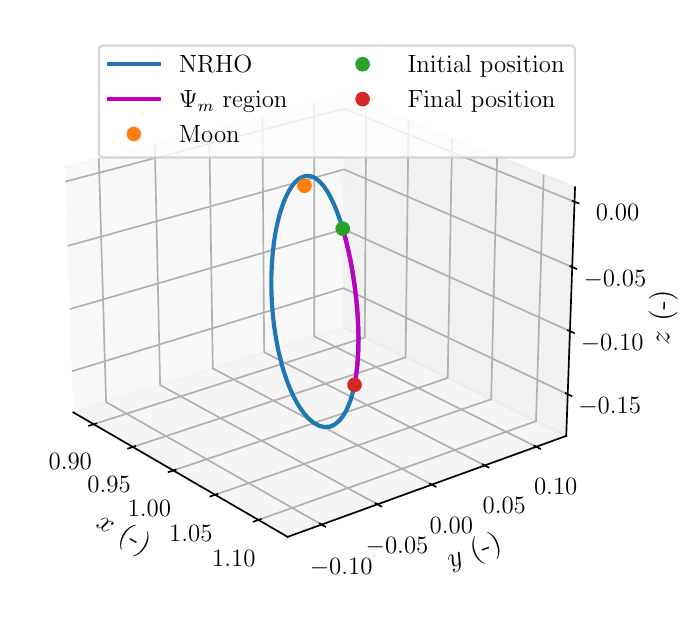}
    \caption{Target's NRHO and time interval of the generated $\Psi_m$ for  rendezvous {scenario} 2.}
    \label{fig:caseA_nrho}
    \end{minipage}
\end{figure}

\subsection{Fixed-Time Optimization Results} \label{sec: results scenatio 1}
The proposed guidance algorithm is applied to solve a large initial separation (about $1500$ \unit{km}) rendezvous case, employing a $4\textsuperscript{th}$ order TSE. The Time of Flight (TOF) is approximately $1.6312$ day, which corresponds to one-quarter of the selected NRHO's period. Figure \ref{fig:caseB_nrho} illustrates the reference trajectory used for the high-order expansion in this example. The chaser's trajectory starts before the pericenter and intersects with the target after this point. For completeness, Table \ref{tab: comp_time_size_stt_daceypy} reports the computational time needed to generate the matrix $\Psi_m$ in the DA framework, along with the estimated RAM usage and the map size in kilobytes when stored as a binary file, averaged over $10$ computations. Note that the map $\Psi_m$ needs to be computed only once on ground and then uploaded onto the spacecraft.

\begin{table}[h]
    \centering 
        \caption{Computational time, RAM, and size for $\Psi_m$ for the scenario $1$.}
    \begin{tabular}{l c c c}
    \toprule
     Order & Time (\unit{s}) & RAM usage (\unit{MB}) & Size (\unit{kB}) \\
    \midrule \midrule     
    $1$ & $17.0821$ & $74.6936$ & $56.8955$ \\
    $2$ & $17.6013$ & $76.6802$ & $234.1093$ \\
    $3$ & $20.2770$ & $83.0992$ & $706.6709$ \\
    $4$ & $31.1817$ & $86.4207$ & $1769.9277$ \\
    \bottomrule
    \end{tabular}
    \label{tab: comp_time_size_stt_daceypy}
\end{table}

\begin{table}[h]
    \centering 
    \caption{Simulation parameters for rendezvous scenario 1.}
    \begin{tabular}{c l}
    \toprule
     \textbf{Parameter} & \textbf{Values}  \\
    \midrule \midrule
    Control interval & $\text{TOF}=1.6312$ days, discretized into $180$ points\\ \midrule
    Relative state at $t_0$ & $\boldsymbol{x}_0 = (1500, -20, 200, -8.9, 13.02, 0)^\top$ (units: $\unit{km}, \unit{m/s}$)\\
    Relative state at $t_f$ & $\boldsymbol{x}_f = (15, 0, 0, 0, 0, 0)^\top$ (units: $\unit{km}, \unit{m/s}$)\\
    \bottomrule
    \end{tabular}
    \label{tab: sim par case D}
\end{table}

\begin{figure}[!h]
    \centering
    \begin{minipage}{0.45\textwidth}
        \centering
        \includegraphics[width=\textwidth]{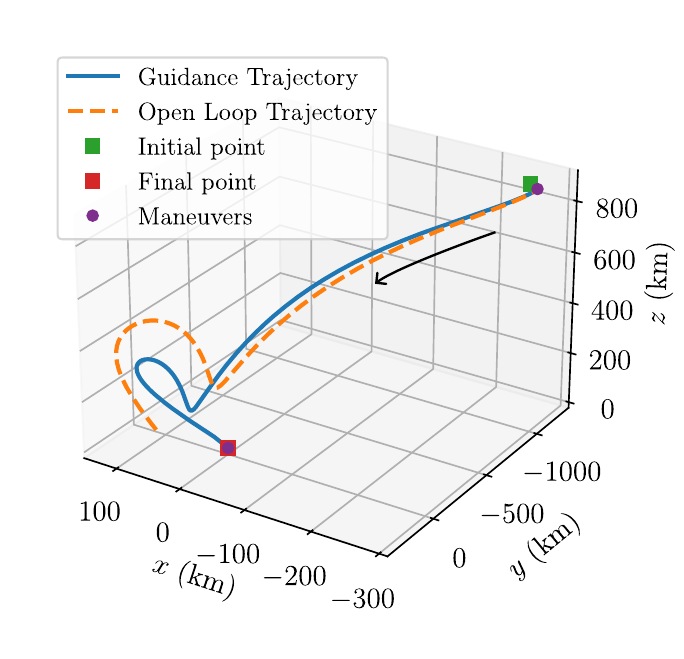}
        \caption{{Linearized relative-motion trajectory in the relative synodic frame for scenario 1.}}
        \label{fig:linea_traj_caseD_lr}
    \end{minipage}
    \hfill
    \begin{minipage}{0.45\textwidth}
        \centering
        \includegraphics[width=\textwidth]{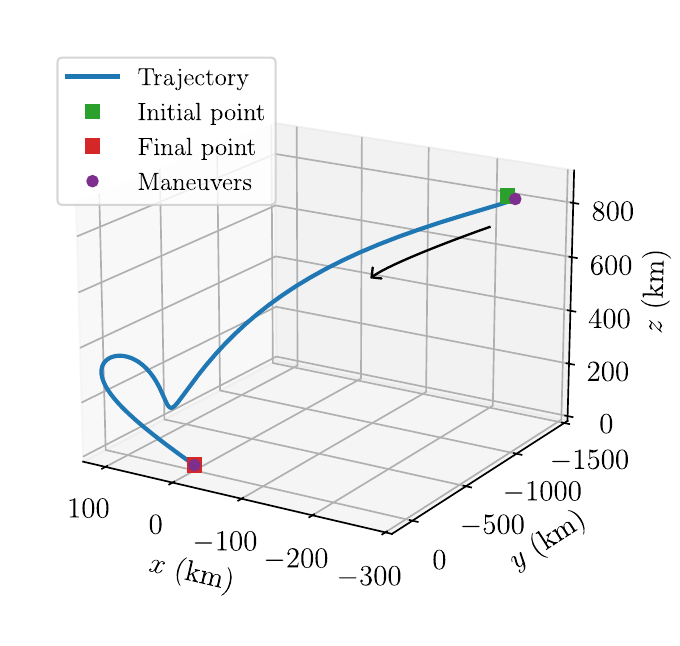}
        \caption{{Relative-motion trajectory in the relative synodic frame for scenario 1.}}
        \label{fig:nonlinea_traj_caseD_lr}
    \end{minipage}
\end{figure}

Table \ref{tab: sim par case D} presents the duration of the TOF, the number of discretized points, and the initial and final states of the chaser, all expressed in the LVLH frame for the CR3BP as defined by \cite{franzini_relative_2019}. The solution obtained by solving the linear version of the problem is illustrated in Figure \ref{fig:linea_traj_caseD_lr}.
The impulsive maneuvers determined by the linear solution cause the trajectory to deviate from the desired path, as can be seen in Figure \ref{fig:linea_traj_caseD_lr}. Therefore, an expansion order greater than $1$ is needed; however, the linear solution can still be used as SCP's initial guess. The SCP guidance is a two-impulse solution with a total $\Delta v$ of $10.74$ \unit{m/s}, and all the other guidance parameters are reported in Table \ref{tab: case D_lr results}. The SCP guidance runtime in the table is the sum of linear guidance runtime and SCP algorithm convergence time.  The chaser's relative trajectory derived from the SCP approach is shown in Figure \ref{fig:nonlinea_traj_caseD_lr}.

\begin{table}[H]
    \centering
    \caption{Guidance results for rendezvous scenario 1.}
    \begin{tabular}{c l}
    \toprule
     \textbf{Parameter} & \textbf{Values}  \\
    \midrule \midrule
       \multirow{2}{*}{Linear guidance}
      & Burn index: $[1, 179]$, $\Delta v = 10.2438$ \unit{m/s},\\ & Runtime = $0.5136$ \unit{s}, Open-loop error $> 225$ km \\ \midrule
    \multirow{3}{*}{SCP guidance} & Burn index: $[1, 179]$, $\Delta v = 10.7495$ \unit{m/s},\\ & Runtime = $6.0355$ \unit{s}, Number of iterations = $4$ \\ & $e_{\text{g}} = (3.1322\e{-7}\; \text{\unit{km}}, 3.4428\e{-10}\; \text{\unit{m/s}}) $ \\ 
    \bottomrule
    \end{tabular}
    \label{tab: case D_lr results}
\end{table}

\begin{figure}[!h]
    \centering
    \begin{minipage}{0.45\textwidth}
        \centering
        \includegraphics[width=\textwidth]{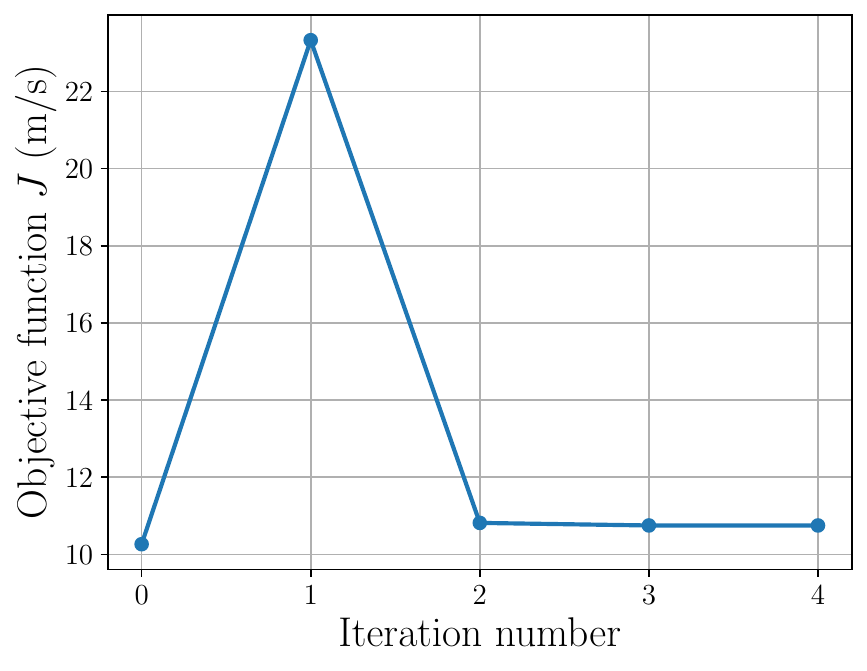}
        \caption{SCP total cost vs. iteration number for scenario 1.}
        \vspace{0.4cm}
        \label{fig:cost_iteration_caseD_lr}
    \end{minipage}
    \hfill
    \begin{minipage}{0.45\textwidth}
        \centering
        \includegraphics[width=\textwidth]{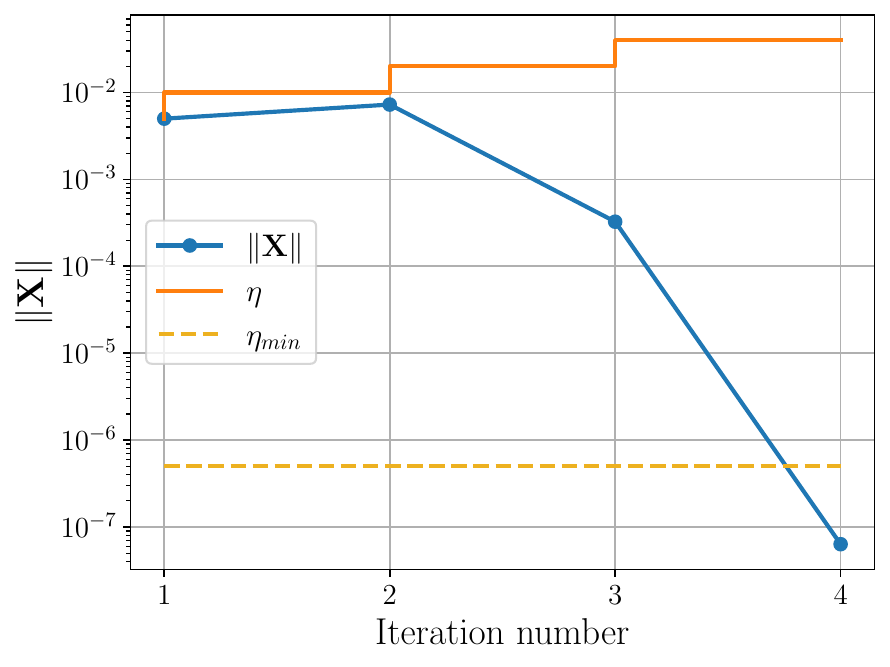}
        \caption{Free variable norm and trust region radius across iteration for scenario 1.}
        \label{fig:norm_caseD_lr}
    \end{minipage}
\end{figure}

\noindent The SCP scheme successfully converges in $4$ iterations, requiring a solve time of $6.04$ \unit{s}. To assess whether a fourth-order expansion for the matrix $\Psi$ is the proper choice, the problem was re-solved under different orders of expansion. Table \ref{tab: order-ol error case Dlr}  reports the open-loop error both for position and velocity. In this rendezvous scenario, the SCP formulation based on a second-order monomial matrix fails to converge to the optimal solution, indicating that an expansion order of at least three is required. An open-loop error lower than  $0.1$ \unit{km} is achieved with an expansion order of $4$. To complete the analysis, Figure \ref{fig:cost_iteration_caseD_lr} and Figure \ref{fig:norm_caseD_lr} report the SCP total cost over iteration and trust region radius update,  indicating stable algorithm performance.

\begin{table}[h]
\centering
\begin{minipage}{0.45\textwidth}
\centering
        \caption{{Open-loop errors for scenario $1$}.}
        \begin{tabular}{c c c}
        \toprule
        \textbf{Order} & \textbf{Position error (\unit{km})} & \textbf{Velocity error (\unit{m/s})} \\
        \midrule \midrule
        $2$ & - & - \\
        $3$ & $6.6310 \e{-1}$ & $6.9586 \e{-3}$ \\
        $4$ & $6.1730 \e{-2}$ & $5.9551 \e{-4}$ \\
        \bottomrule
        \end{tabular}
        \label{tab: order-ol error case Dlr}
\end{minipage}
\hfill
\begin{minipage}{0.45\textwidth}
\centering
        \caption{{Open-loop errors for scenario $2$}.}
        \begin{tabular}{c c c}
        \toprule
        \textbf{Order} & \textbf{Position error (\unit{km})} & \textbf{Velocity error (\unit{m/s})} \\       
        \midrule \midrule
        $2$ & $5.1044 \e{-2}$ & $4.2447 \e{-4}$ \\
        $3$ & $7.5950 \e{-3}$ & $5.4142 \e{-5}$ \\
        $4$ & $2.9924 \e{-3}$ & $1.8952 \e{-6}$ \\
        \bottomrule
        \end{tabular}
        \label{tab: order-ol error case A}
\end{minipage}
\end{table}

\begin{figure}[!h]
    \centering
    \begin{minipage}{0.45\textwidth}
        \centering
        \includegraphics[width=\textwidth]{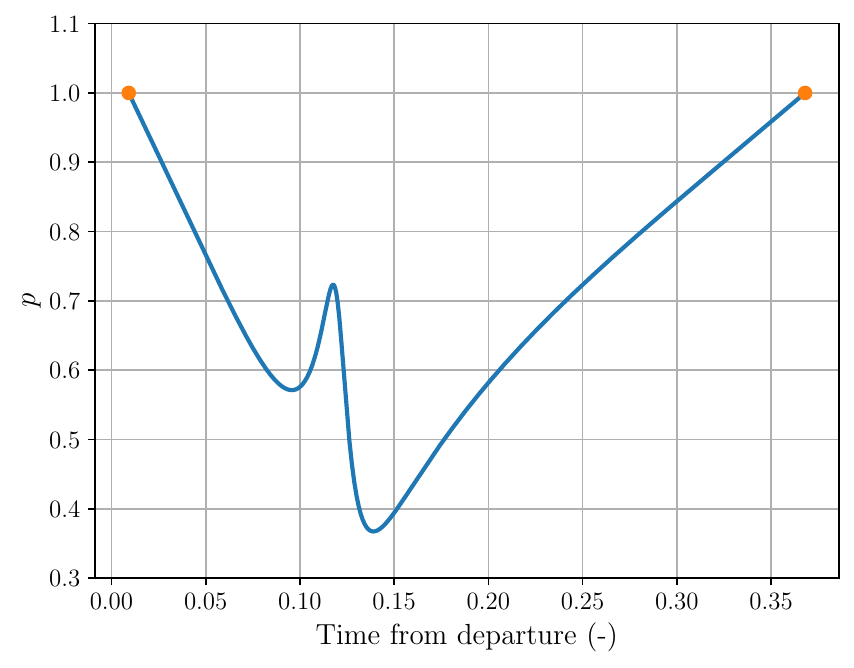}
        \caption{Primer vector norm along the trajectory for scenario 1.}
        \vspace{0.4cm}
        \label{fig:primer-vector-norm case Dlr}
    \end{minipage}
    \hfill
    \begin{minipage}{0.45\textwidth}
        \centering
        \includegraphics[width=\textwidth]{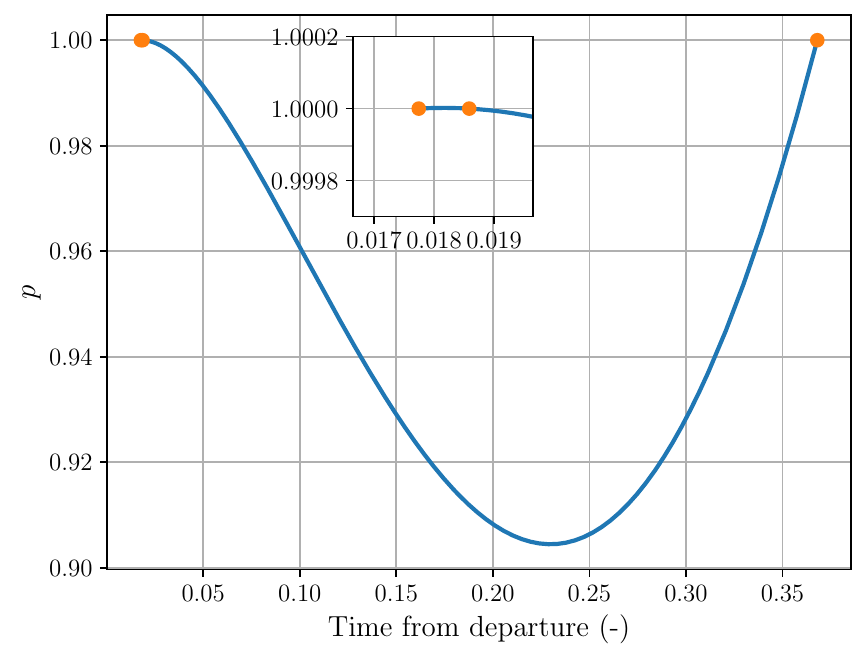}
        \caption{Primer vector norm along the trajectory for scenario 3.}
        \label{fig:primer-vector-norm case A unconstrained}
    \end{minipage}
\end{figure}

The optimality of the guidance solution generated by the proposed monomial-based method is assessed \textit{a posteriori} using primer-vector theory \cite{primer-vector-1963}. For a solution to be optimal, the primer vector must satisfy the following conditions: (i) its magnitude must be less or equal to unity at all times, (ii) it must attain unit magnitude at each impulse time, and (iii) the direction of each impulse velocity increment must be aligned with the primer vector at the corresponding maneuver time. 
In this work, the state and the costate histories required to construct the primer vector are obtained from the guidance solution computed using the monomial-based framework. Once the sequence of impulsive maneuvers and their corresponding time are determined, the associated costate dynamics are reconstructed by propagating the adjoint equations between impulses \cite{bucchinio-primer-vector}. 
Figure \ref{fig:primer-vector-norm case Dlr} shows the time history of the primer vector norm corresponding to the computed guidance solution. As observed, the norm remains strictly below unity between maneuvers. This behavior is fully consistent with the conditions of optimality described by primer vector theory. 

\subsection{State Constraints in Nonlinear Convex Rendezvous Guidance} \label{sec:Stt constra}
Figure \ref{fig:caseA_nrho} illustrates the portion of the target's NRHO used as a reference to generate the monomial matrix  $\Psi_m$ for this second {scenario}. The chaser starts with an offset from the target of about $62$ \unit{km} on a point after the pericenter and reaches the target prior to the orbit's apocenter. Table \ref{tab: comp_time_size_stt_daceypy_case_cons} reports the computational time to generate the monomial matrix, RAM usage, and storage size for the current case for different expansion orders.  As in the previous example, the TOF is fixed and it is unchanged. A time-dependent minimum-range constraint, represented as a step function, is applied to the chaser trajectory.  The minimum range is chosen to avoid drastic trajectory changes, encouraging smoother adjustments instead. The approach cone has a semi-aperture angle of $\gamma = \ang{15}$, aligned with the V-bar in the LVLH frame, and is active only during the last $20\,\%$ of the trajectory nodes.

\begin{table}[h]
    \centering 
    \caption{Computational time, RAM, and size for $\Psi_m$ for the  scenario 2.}
    \begin{tabular}{l c c c}
    \toprule
     Order & Time (\unit{s}) & RAM usage (\unit{MB}) & Size (\unit{kB}) \\
    \midrule \midrule     
    $1$ & $12.8830$ & $74.3523$ & $31.6611$ \\
    $2$ & $13.4391$ & $76.0769$ & $130.1074$ \\
    $3$ & $15.5433$ & $82.1976$ & $392.6426$ \\
    $4$ & $24.0456$ & $83.9932$ & $983.3467$ \\
    \bottomrule
    \end{tabular}
    \label{tab: comp_time_size_stt_daceypy_case_cons}
\end{table}
The minimum-range and approach cone constraints are represented as a $4\textsuperscript{th}$ order expansion, as formulated in Section \ref{sec:scp cons}. The minimum-range constraint map takes $23.4799$ \unit{s} to be generated and it occupies $133.4681$ \unit{kB}, while for the approach cone, $24.0365$ \unit{s} are needed with the size of $33.4580$ \unit{kB}.

\begin{figure}[H]
\includegraphics[width=.49\linewidth]{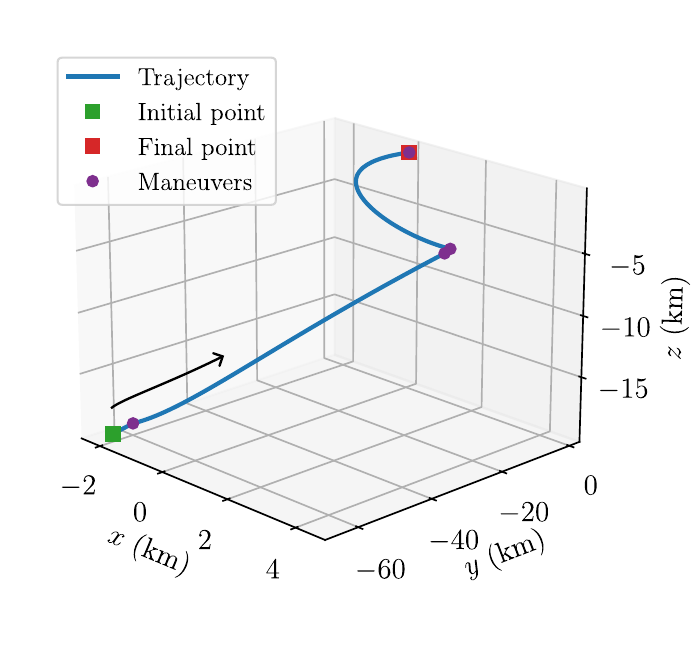}\hfill
\includegraphics[width=.49\linewidth]{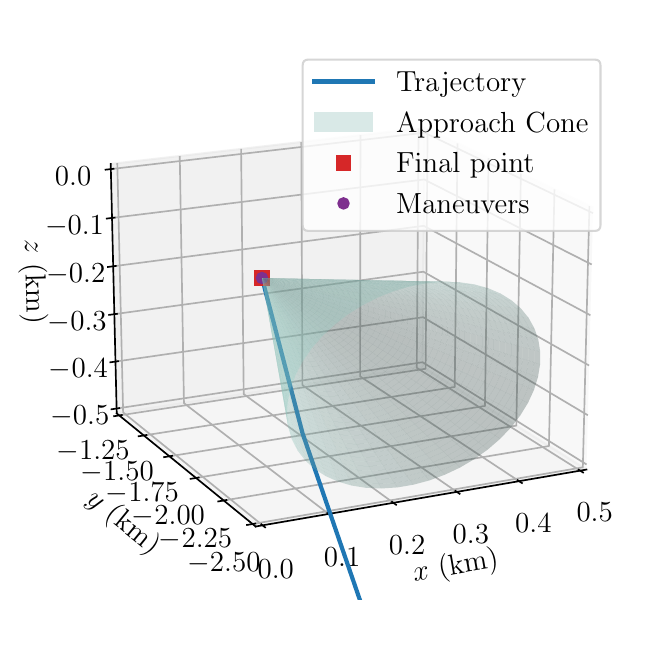}
        \caption{Relative-motion trajectory in the relative synodic frame for scenario 2.}
\label{fig:nonlinea_traj_caseA_cons}
\end{figure}

\begin{figure}[!h]
    \centering
    \begin{minipage}{0.45\textwidth}
        \centering
        \includegraphics[width=\textwidth]{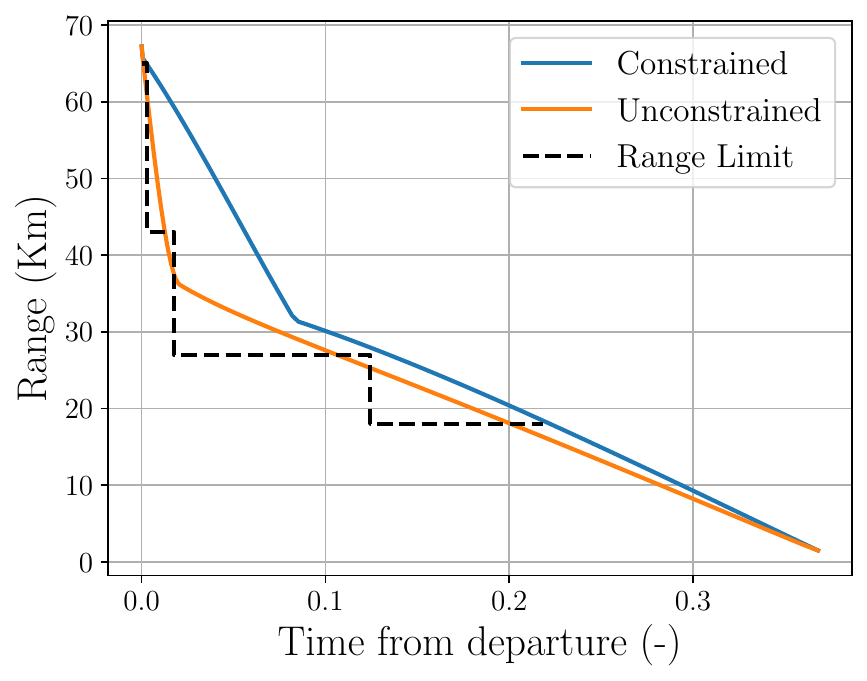}
        \caption{Inter-spacecraft range for the rendezvous {scenario} 2.}
        \label{fig:range}
    \end{minipage}
    \hfill
    \begin{minipage}{0.45\textwidth}
        \centering
        \includegraphics[width=\textwidth]{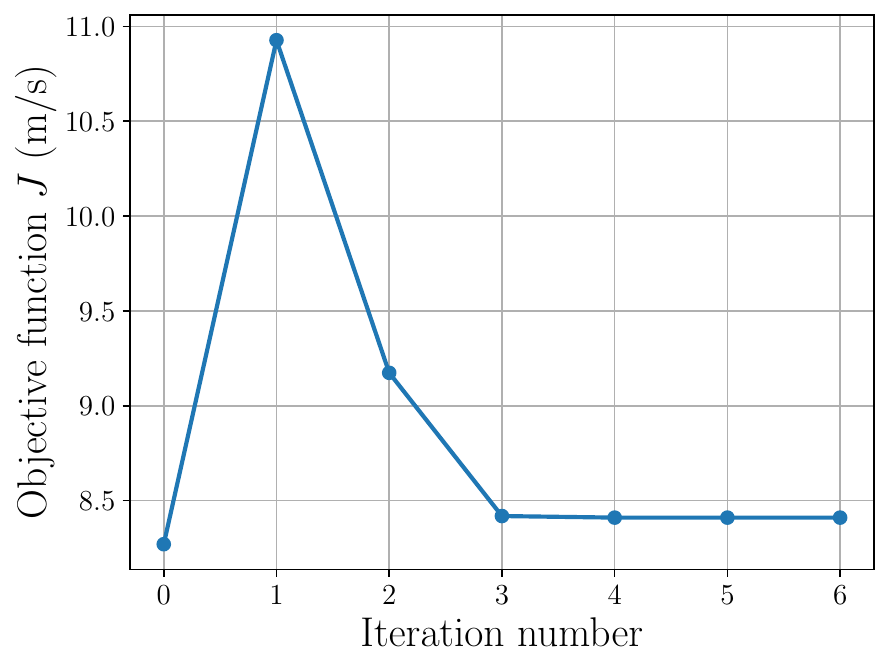}
        \caption{SCP total cost vs. iteration number for the rendezvous  {scenario} 2.}
        \label{fig:cost_cons_caseA}
    \end{minipage}
\end{figure}

The simulation parameters are reported in Table \ref{tab: sim par case D} while the SCP parameters, Table \ref{tab: SCP parameter}, are unchanged compared to the rendezvous scenario $1$, but we re-tune $w=2500$ to improve convergence characteristics. Slack variables introduced by constraints linearization are equally penalized in the objective function with a penalty weight of $7000$. The SCP algorithm converges in $6$ iterations, with the cost decreasing across iterations as shown in Figure \ref{fig:cost_cons_caseA}. 
The SCP-constrained guidance solution consists of $4$ impulsive maneuvers, with a total $\Delta v$ of $8.40$ \unit{m/s}. The SCP results are summarized in Table \ref{tab: case A cons results} {, while the open-loop errors obtained for different expansion order are reported in Table \ref{tab: order-ol error case A}. In contrast to the scenario presented in Section \ref{sec: results scenatio 1}, which involved a much larger initial separation, an expansion order of two is sufficient in this case to achieve a final open-loop error below $0.1$ \unit{km}}. Figure \ref{fig:nonlinea_traj_caseA_cons} reports the trajectory resulting from the constrained nonlinear convex guidance formulated with the monomial parameterization for the minimum range and the approach cone. {On the right of the figure, a zoomed-in view of the trajectory's endpoint highlights its position within the approach cone.} The step-function minimum range and inter-spacecraft range over time are illustrated in Figure \ref{fig:range}.

\begin{table}[h]
    \centering
    \caption{Simulation parameters for rendezvous scenario 2.}
    \begin{tabular}{c l}
    \toprule
     \textbf{Parameter} & \textbf{Values}  \\
    \midrule \midrule
    Control interval & $\text{TOF}=1.6312$ days, discretized into $100$ points\\ \midrule
    Relative state at $t_0$ & $\boldsymbol{x}_0 = (62, -7, 25, -6.59, 3.46, 0)^\top$ (units: $\unit{km}, \unit{m/s}$)\\
    Relative state at $t_f$ & $\boldsymbol{x}_f = (1.5, 0, 0, 0, 0, 0)^\top$ (units: $\unit{km}, \unit{m/s}$)\\ 
    \bottomrule
    \end{tabular}
    \label{tab: sim para case A}
\end{table}

\begin{table}[h]
    \centering
    \caption{Guidance results for rendezvous scenario 2.}
    \begin{tabular}{c l}
    \toprule
     \textbf{Parameter} & \textbf{Values}  \\
    \midrule \midrule
       \multirow{2}{*}{Linear guidance}
      & Burn index: $[15, 16, 99]$, $\Delta v = 8.2703$ \unit{m/s},\\ & Runtime = $0.2299$ \unit{s}  \\ \midrule
    \multirow{3}{*}{SCP guidance} & Burn index: $[1, 47, 48, 99]$, $\Delta v = 8.3985$ \unit{m/s},\\ & Runtime = $4.8404$ \unit{s}, Number of iterations = $6$ \\ & $e_{\text{g}} = (6.0703\e{-6}\; \text{\unit{km}}, 3.1587\e{-9}\; \text{\unit{m/s}}) $ \\ 
    \bottomrule
    \end{tabular}
    \label{tab: case A cons results}
\end{table}

\subsection{Free-Time Optimization Results} \label{sec: free time results}
Now we provide a simple demonstration of the free-time rendezvous. The outer-bound time range, discretized in $100$ nodes, is set as the TOF used in the previous cases for simplicity, and the map used in Section \ref{sec:Stt constra} is interpolated with a natural spline to provide a smooth approximation of  $\Psi_m(t_0, t)$ and its derivatives over the time range.  
Since a limited number of impulses characterizes the fixed-time rendezvous solutions, we set $N_l=6$ to limit problem dimensionality, and to test the free-time approach's ability to adjust burn nodes, we start with an initial SCP guess with a TOF of $0.1406$ \unit{days}, discretized over seven nodes (including the initial node) indexed as $[0, 1, 5, 10, 14, 19, 24]$ in $[t_0, t_\text{max}]$. SCP parameters for the free-time case match those in Table \ref{tab: SCP parameter}. 

\begin{table}[h]
    \centering
    \caption{Guidance Results for rendezvous scenario 3. }
    \begin{tabular}{ c l }
    \toprule
     \textbf{Parameter} & \textbf{Values}  \\
    \midrule \midrule
    \multirow{2}{*}{Linear guidance}
      & $\text{TOF}=0.1406$ days, discretized into $7$ points,\\ & Runtime = $0.0197$ \unit{s}, $\Delta v = 10.0144 \unit{m/s}$   \\ \midrule
    \multirow{3}{*}{SCP guidance} & $\text{TOF}=1.6312$ days, discretized into $7$ points,\\ & Near index: $[0, 2, 4, 15, 32, 98, 99]$, $\Delta v = 8.2709$ \unit{m/s},\\ & Runtime = $3.9725$ \unit{s}, Number of iterations = $20$ \\ \hdashline 
    \multirow{4}{*}{SCP correction} & Burn index: $[15, 99]$, $\Delta v = 8.2711$ \unit{m/s},\\ & Runtime = $0.0865$ \unit{s}, Number of iterations = $2$ \\ & $e_{\text{g}} = (8.1715\e{-5}\; \text{\unit{km}}, 5.0873\e{-8}\; \text{\unit{m/s}}) $ \\ & $e_{\text{op}}=(2.4815\e{-2}\; \text{\unit{km}}, 1.9692\e{-5}\; \text{\unit{m/s}})$\\ \midrule
     \multirow{3}{*}{SCP fixed-final time} & $\text{TOF}=1.6312$ days, discretized into $100$ points,\\ & Burn index: $[15, 16, 99]$, $\Delta v = 8.2709$ \unit{m/s},\\ & Runtime = $3.1792$ \unit{s}, Number of iterations = $4$ \\ 
     \bottomrule
    \end{tabular}
    \label{tab: case A free results}
\end{table}

\begin{figure}[h]
    \centering
    \begin{minipage}{0.45\textwidth}
        \centering
        \includegraphics[width=\textwidth]{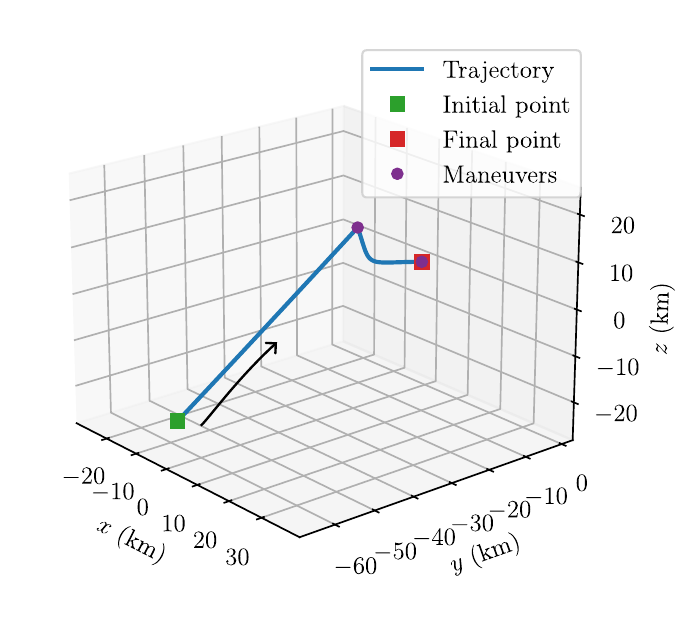}
        \caption{{Relative-motion trajectory in the relative synodic frame for rendezvous scenario 3.}}
        \label{fig:case free traj}
    \end{minipage}
    \hfill
    \begin{minipage}{0.45\textwidth}
        \centering
        \includegraphics[width=\textwidth]{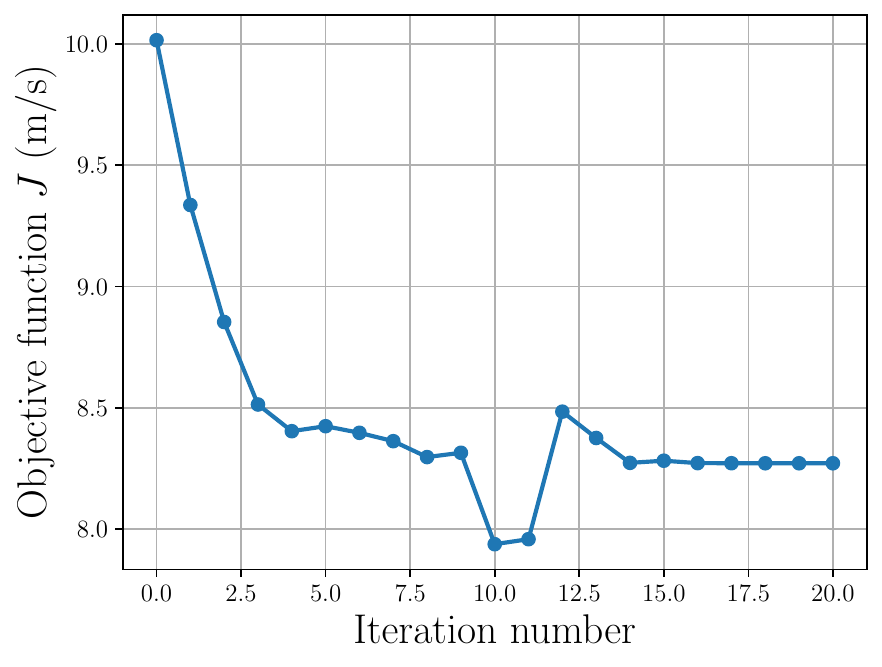}
        \caption{SCP total cost vs. iteration number for rendezvous scenario 3.}
        \label{fig:caseA free cost}
    \end{minipage}
\end{figure}

The case analyzed in the context of free-time SCP has an initial separation of approximately $62$ \unit{km}, and the solution can extend over the entire purple region of Figure \ref{fig:caseA_nrho}. The simulation parameters for the initial and final states of the chaser are provided in Table \ref{tab: sim para case A}. The solution of the linear approach with a starting TOF of $0.1406$ \unit{days} acts as the initial guess for the free-time SCP. Table \ref{tab: case A free results} presents the solution of the converged free-time SCP and the post-processed trajectory using the fixed-time SCP discussed in earlier sections. Both guidance and open-loop errors remain within the desired value. Although the initial TOF is small, the optimal solution for the free-time SCP uses the entire available time range, resulting in a TOF of $1.6312$ days. This behavior has been consistently observed across all tested free-time rendezvous scenarios, which have not been reported in this paper. In fact, for the same initial separation, a longer TOF allows for smaller and more efficient burns, which results in lower fuel consumption due to the enhanced exploitation of natural dynamics between two maneuvers. This solution has then been compared with the results of the fixed-final time approach over the same NRHO's region resulting from the free-time approach. The two strategies require a similar $\Delta v$. {As in rendezvous scenario $1$ (Section \ref{sec: results scenatio 1}), the optimality of the solution is verified \textit{a posteriori} using primer-vector theory. Figure \ref{fig:primer-vector-norm case A unconstrained} shows that, for the three-impulse solution, the norm of the primer vector remains less than or equal to unity over the entire trajectory.} The trajectory for the free-time guidance is shown in Figure \ref{fig:case free traj} and the corresponding cost for the free-time SCP is highlighted in Figure  \ref{fig:caseA free cost}, displaying the cost reduction over SCP iterations.

Analyzing the free-time guidance results in Table \ref{tab: case A free results}; as expected, the postprocess solution exhibits a slightly higher $\Delta v$ than the optimal free-time SCP since the correction yields a {slightly} suboptimal trajectory. This increment is the natural consequence of enforcing burn times at discrete nodes. Moreover, the free-time guidance requires approximately $16\,\%$ more computation time of the corresponding fixed-time case. In this scenario, the algorithm not only determines the impulsive maneuvers but also modifies the time of flight. As expected, using fewer nodes ($N_l \ll N$), the linear guidance results in significantly lower computation time compared to a linear guidance solution with $100$ nodes. While this analysis prioritizes a solver that freely adjusts the TOF, a larger initial TOF, comparable to the $\Psi_m$ region, would yield a faster convergence in the free-time case than the fixed-time approach.

\subsection{Reducing Computational Overhead for the Monomial Method and Comparison with Classical SCP Approach}
\begin{figure}[h]
    \centering
    \includegraphics[width=1\linewidth]{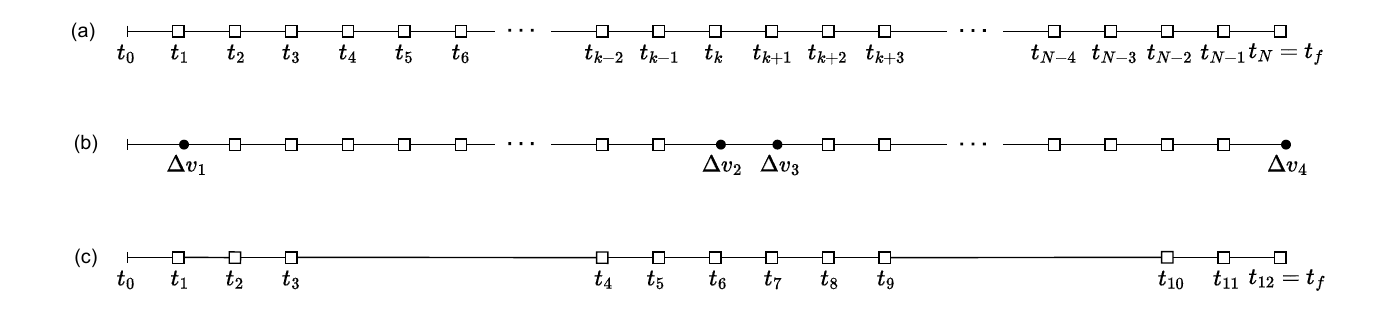}
    \caption{Illustration of reducing time-grid nodes ((a) $\boldsymbol{\rightarrow}$ (b) $\boldsymbol{\rightarrow}$ (c)) for a $\boldsymbol{4}$-impulse linear guess.}
    \label{fig:time_grid_reduced}
\end{figure}

\begin{figure}
    \centering
    \includegraphics[width=0.5\linewidth]{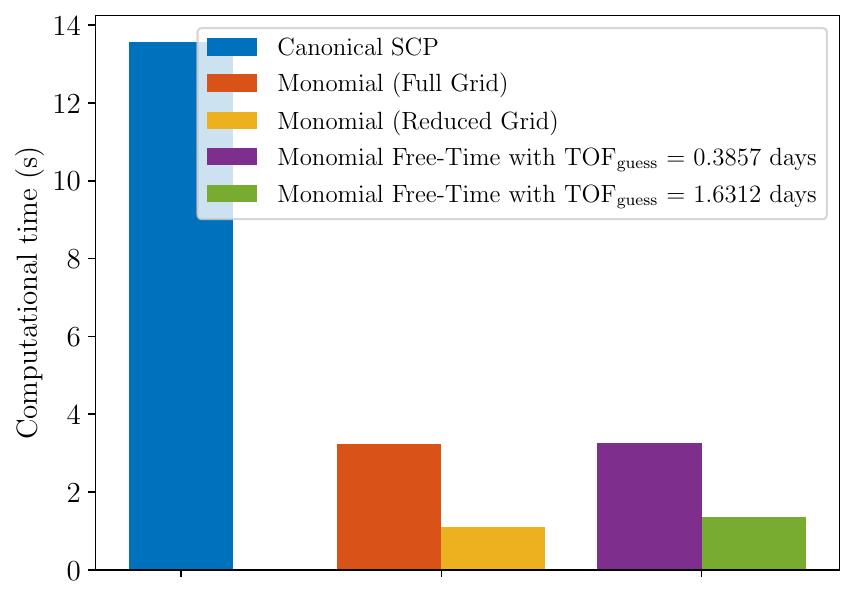}
    \caption{Comparison of computational time for different guidance methods.}
    \label{fig:histogram comparison}
\end{figure}

To reduce the computational time of the fixed-time monomial method, the solution from the initial guess of the SCP is used to decrease the number of nodes in the time grid. This is because we observe that the optimal maneuver times in the nonlinear solution do not change much from the linear solution. This process refines the discretization by focusing on the most relevant trajectory points. The key idea is to retain only a user-defined number of nodes adjacent to the impulsive nodes obtained from the linear solution. Figure \ref{fig:time_grid_reduced} illustrates this procedure for a four-impulse initial guess, where only two nodes before and after each impulse node are selected. By applying this approach, the total number of nodes, and consequently the computational time, is reduced, as shown in Figure \ref{fig:histogram comparison}. This method remains effective and converges to the same solution as the full-grid approach, as long as the linear solution produces maneuvering nodes close to those obtained by the nonlinear algorithm. 

In Section \ref{sec: free time results}, we applied the free-time approach with a very small initial TOF guess to demonstrate the algorithm’s feasibility. However, in this comparison, we increased the initial TOF guess to evaluate whether the free-time approach could outperform or match the fixed-time method in terms of efficiency. We tested two different initial TOF values: one set to $0.3857$ days, which converged in a similar time to the monomial method with $100$ nodes, and another using the maximum admissible TOF, corresponding to the NRHO discretization region length. In the latter case, the algorithm converged in a comparable time to the monomial method with a reduced time grid. However, unlike the time grid reduction approach (whose solution accuracy depends on the quality of the linear solution) the free-time approach allows for maneuvering nodes to shift, starting from any initial guess.

Finally, the performance of the suggested monomial parameterization is compared with a standard SCP method for impulsive trajectory design in order to demonstrate its computational and practical advantages. The canonical SCP method follows the approach outlined in \cite{malyuta_convex_2022}, applied to an impulsive rendezvous problem. This method employs STMs, computed along the reference trajectory using the variational equation, to approximate the state deviation of the spacecraft from the reference and formulate a convex subproblem that is solved iteratively \cite{impulsiveSTMSCP_rizza}. The time domain is discretized using $N$ nodes, consistent with the monomial-based method, and the SCP is initialized using a linear interpolation between the initial and final states. A key computational burden of the canonical SCP approach arises from the need to repeatedly integrate the dynamics and the associated variational equations between all nodes. 
Each integrator step requires multiple evaluations of the nonlinear dynamics, with the exact number depending on the integrator’s order and tolerances set, and this computational load is incurred at each iteration of SCP. The monomial parameterization method, by contrast, offers a more efficient alternative, though it still requires nonlinear operations between SCP iterations, as it updates the monomial basis at each trajectory node. The monomial approach operates with a limited number of discretization nodes, typically ranging from a few to around a hundred. This reduces the number of nonlinear operations by shifting the computational burden from evaluating nonlinear dynamics at each internal integrator's node to evaluating monomial bases at each parameterization node. As shown in Figure \ref{fig:histogram comparison}, all the monomial-based approaches presented in this paper require significantly less computational time than the canonical SCP approach. Moreover, Table \ref{tab: deltaV_comparison} shows that all the considered guidance strategies achieve essentially identical values of total $\Delta v$, demonstrating that the proposed monomial parameterization preserves the solution optimality while significantly reducing the computational time. 

\begin{table}[h]
    \centering 
    \caption{Comparison of total $\boldsymbol{\Delta v}$ for different guidance methods.}
    \begin{tabular}{l c}
    \toprule
    Guidance method & $\Delta v$ (\unit{m/s}) \\
    \midrule \midrule     
     Canonical SCP & $8.2710$ \\
     Monomial (Full Grid) & $8.2755$ \\
     Monomial (Reduced Grid) & $8.2703$ \\
     Monomial Free Time with $\text{TOF}_\text{guess} = 0.3857$ days & $8.2711$ \\
     Monomial Free Time with $\text{TOF}_\text{guess} = 1.6312$ days & $8.2712$ \\
    \bottomrule
    \end{tabular}
    \label{tab: deltaV_comparison}
\end{table}


\section{Conclusion} \label{sec: conclusion}
A nonlinear convex guidance method based on high-order expansion of dynamics is proposed for addressing fixed-final time trajectory generation problem near a reference. Using monomial parameterization, both dynamics and constraints are approximated, and the problem is posed on elements of a non-Euclidean sub-manifold.  This approach allows for trajectory optimization formulated as a path-planning problem with a low real-time computational load, making it ideal for on-board implementation with limited computational resources. Indeed, thanks to monomial parameterization, instead of having to repeatedly evaluate the nonlinear dynamics, as in the case of the classical SCP method, it is sufficient to update the monomial equations between iterations, which results in a much smaller set of numerical computations. To generalize this method to free-final time problems, burn times have been included in the optimization variable. Here, burn times can be chosen arbitrarily over the time window but are later enforced in post-processing to be associated with the discrete nodes, thus removing interpolation errors while retaining only truncation ones. Additionally, extending to a free-final time reduces the number of optimization variables, potentially saving computational time. The method is particularly well-suited for spacecraft relative motion and rendezvous problems, where the reference trajectory, needed to generate the matrix $\Psi_m$, is the orbit of the target. However, the approach is general and can be applied to any dynamics, even ephemeris-like. 

The methodology is tested in rendezvous problems within the framework of the cislunar domain with the target orbiting on an NRHO. The nonlinear convex guidance proves to be effective in different rendezvous scenarios, with and without a fixed final time, that differ in terms of nonlinearity level and different ranges of initial separation from $62$ \unit{km} to $1500$ \unit{km}. The SCP-based guidance takes a few seconds to run with a longer solving time for the higher nonlinear case, and it is able to generate a fuel-optimal solution that reaches the final state precisely with small guidance and open-loop errors. For the unconstrained cases, the optimality of the computed solutions has been further verified \textit{a posteriori} using primer vector theory. Finally, the approach demonstrates practical onboard feasibility, as a Python prototype executes within seconds on a modern computing platform.

\appendix
\section{Monomial Parameterization and STTs Equivalence} \label{app:stt and mon matrix equi}
\begin{figure}[h]
    \centering
    \includegraphics[width=\textwidth]{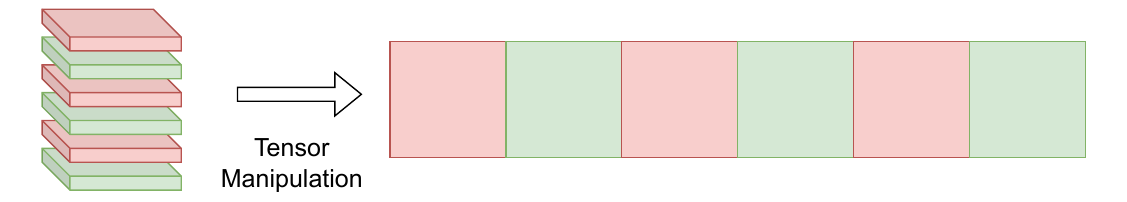}
    \caption{Example of horizontal slices for a $\mathbf{3}$rd-order tensor with its manipulation. }
    \label{fig:Tensor slice manipulation}
\end{figure}

This appendix shows that through manipulation, STTs directly lead to monomial parameterization; hence, the two approaches are equivalent. There are two main differences between them. Firstly, the matrix $\Psi_m$ has dimension of $n\times K_m$, while the STTs consist of a collection of $m$ tensors of dimension $n^2, n^3, \dots, n^{m+1}$. However, a tensor manipulation procedure can reformulate them into a single matrix with $n$ rows.  Given a tensor of generic order, a canonical tensor manipulation is proposed: all the 'slices' (or sub-matrices) of the tensor are stacked side by side to form a matrix. A 'slice' is a two-dimensional section of a tensor, defined by fixing all but two indices \cite{kolda_tensor_2009}. A general $l$-th order tensor $\Lambda^{(l)}$ can be represent as a two dimensional matrix of size $n \times n^{l-1}$, assuming that the length of each dimension of $\Lambda^{(l)}$ is  equal to $n$, using the following relationship \cite{spreen_automated_2017}:
    \begin{equation}
        \Psi^{(l)}(i_1, i_2 + n(i_3-1) + n^2(i_4-1) + \dots + n^{l-2}(i_l - 1) = \Lambda^{(l)}(i_1, i_2, \dots , i_l)
    \end{equation}
where the index $i_1, i_2, \dots, i_l$ span, separately, from $1$ to $n$, and the superscript $(l)$ denotes the order of the tensor. An example of the "horizontal" slices for a third-order tensor, along with the tensor manipulation in which each slice is stacked side by side, is reported in Figure \ref{fig:Tensor slice manipulation}.

\begin{figure}[h]
    \centering
    \includegraphics[width=0.8\linewidth]{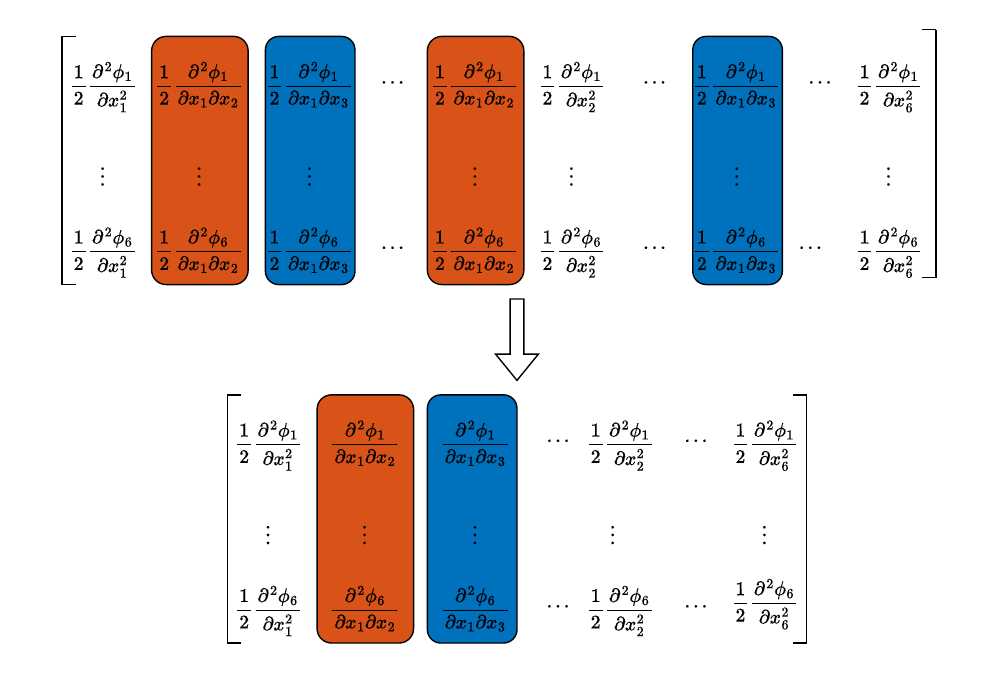}
    \caption{Example of combining redundant terms for a manipulated second order STT. }
    \label{fig:stt2mani}
\end{figure}

\noindent Thanks to this tensor manipulation, given the first $m$ STTs of a system of dimension $n$ and joining each manipulated tensor side by side, it is possible to construct a $n \times (n + n^2 + \dots n^{m-1})$ matrix. Even after the manipulation of the STTs into a matrix, it still differs from the matrix $\Psi_m$. The difference lies in the number of columns because the STT collects redundant terms. In the STT, mixed derivatives are equal (e.g. $\partial \phi_j / (\partial x_1 \partial x_2) = \partial \phi_j / (\partial x_2 \partial x_1)$ $\forall j=1, \dots, n$). Hence, by combining redundant rows in the manipulated STTs, the consequence is a $n\times K_m$ matrix that is exactly the matrix $\Psi_m$, which is the minimal representation for the approximation written in a quasi-linear form. Figure \ref{fig:stt2mani} illustrates an example of the process of combining redundant terms for a manipulated second-order STT for a $6$-dimensional system. Starting from $36$ cols, the manipulated STT reduces to $21$ cols. Contrary to the original definition of STT \cite{park_nonlinear_2006}, the coefficients that normalize the derivatives multiply each element of the manipulated STT.

\bibliography{references}

@misc{pyaudi-repository,
  author    = {Izzo, Dario and Biscani, Francesco},
  title     = {{Audi} and {Pyaudi} documentation},
  year      = {2018},
  url       = {https://darioizzo.github.io/audi/index.html},
  note      = {Accessed: 2025-04-15}
}

@misc{daceypy-repository,
  author    = {Purpura, Giovanni and Maestrini, Michele},
  title     = {{DACEyPy} documentation},
  year      = {2018},
  url       = {https://github.com/giovannipurpura/daceypy},
  note      = {Accessed: 2025-04-15}
}

@incollection{cellier_continuous_2006,
	title = {Continuous {System} {Simulation}},
    isbn = {978-0-387-30260-7},
	language = {en},
	urldate = {2024-12-30},
    publisher={Springer New York, NY},
    author={Cellier, Fran{\c{c}}ois E and Kofman, Ernesto},
    year = {2006},
    pages = {25--164},
	address = {Boston, MA},
    doi = {10.1007/0-387-30260-3}
}

@article{diamond_cvxpy_2016,
	title = {{CVXPY}: {A} {Python}-{Embedded} {Modeling} {Language} for {Convex} {Optimization}},
	volume = {17},
	url = {https://stanford.edu/~boyd/papers/cvxpy_paper},
	language = {en},
	number = {83},
	journal = {Journal of Machine Learning Research},
	author = {Diamond, Steven and Boyd, Stephen},
	year = {2016},
	pages = {1--5},
}

@incollection{szebehely_chapter_1967,
	title = {Chapter 1 - {Description} of the {Restricted} {Problem}},
	isbn = {978-0-12-395732-0},
	url = {https://linkinghub.elsevier.com/retrieve/pii/B9780123957320X50016},
	language = {en},
	urldate = {2024-12-30},
	booktitle = {Theory of {Orbit}},
	publisher = {Elsevier},
	author = {Szebehely, Victor},
	year = {1967},
}

@article{di_lizia_high_2014-1,
	title = {High order optimal feedback control of space trajectories with bounded control},
	volume = {94},
	issn = {0094-5765},
	url = {https://www.sciencedirect.com/science/article/pii/S0094576513000593},
	doi = {10.1016/j.actaastro.2013.02.011},
	number = {1},
	urldate = {2025-04-08},
	journal = {Acta Astronautica},
	author = {Di Lizia, P. and Armellin, R. and Morselli, A. and Bernelli-Zazzera, F.},
	month = jan,
	year = {2014},
	pages = {383--394},
}

@incollection{berz_chapter_1999,
	address = {San Diego, CA},
	series = {Advances in imaging and electron physics},
	title = {Chapter 2—{Differential} {Algebraic} {Techniques}},
	isbn = {978-0-12-014750-2},
	language = {eng},
	number = {108},
	booktitle = {Modern map methods in particle beam physics},
	publisher = {Academic Press},
	author = {Berz, Martin},
	year = {1999},
    pages = {81--117}
}

@book{berkovitz_optimal_1974,
	address = {New York, NY},
	series = {Applied {Mathematical} {Sciences}},
	title = {Optimal {Control} {Theory}},
	volume = {12},
	copyright = {http://www.springer.com/tdm},
	isbn = {978-1-4419-2804-7 978-1-4757-6097-2},
	url = {http://link.springer.com/10.1007/978-1-4757-6097-2},
	urldate = {2025-02-06},
	publisher = {Springer},
	author = {Berkovitz, Leonard D.},
	editor = {John, Fritz and Sirovich, Lawrence and LaSalle, Joseph P. and Whitham, Gerald B.},
	year = {1974},
}

@article{boone_orbital_2021,
	title = {Orbital {Guidance} {Using} {Higher}-{Order} {State} {Transition} {Tensors}},
	volume = {44},
	issn = {0731-5090},
	url = {https://arc.aiaa.org/doi/10.2514/1.G005493},
	doi = {10.2514/1.G005493},
	number = {3},
	urldate = {2024-02-29},
	journal = {Journal of Guidance, Control, and Dynamics},
	author = {Boone, Spencer and McMahon, Jay},
	month = mar,
	year = {2021},
	pages = {493--504},
}

@inproceedings{serra-primer-vector,
  TITLE = {{Fuel-optimal impulsive fixed-time trajectories in the linearized circular restricted 3-body-problem}},
  AUTHOR = {Serra, Romain and Arzelier, Denis and Br{\'e}hard, Florent and Joldes, Mioara},
  URL = {https://hal.science/hal-01830253},
  BOOKTITLE = {{CSA/IAF Special issue IAF Astrodynamics Symposium (69TH international astronautical congress)}},
  ADDRESS = {Bremen, Germany},
  ORGANIZATION = {{International Astronautical Federation}},
  HAL_LOCAL_REFERENCE = {Rapport LAAS n{\textdegree} 18214},
  SERIES = {CSA/IAF Special issue IAF Astrodynamics Symposium (69TH international astronautical congress)},
  PAGES = {1-9},
  YEAR = {2018},
  MONTH = Oct,
  HAL_ID = {hal-01830253},
  HAL_VERSION = {v1},
}

@article{bucchinio-primer-vector,
title = {Optimal time-fixed impulsive non-Keplerian orbit to orbit transfer algorithm based on primer vector theory},
journal = {Communications in Nonlinear Science and Numerical Simulation},
volume = {124},
pages = {107307},
year = {2023},
issn = {1007-5704},
doi = {https://doi.org/10.1016/j.cnsns.2023.107307},
url = {https://www.sciencedirect.com/science/article/pii/S1007570423002253},
author = {Giordana Bucchioni and Gabriele Gemignani and Francesco Lombardi and Andrea Bellome and João Pedro Fernandes Leitão and Stèphanie Lizy-Destrez and Mario Innocenti}
}

@misc{goulart_clarabel_2024,
	title = {Clarabel: {An} interior-point solver for conic programs with quadratic objectives},
	shorttitle = {Clarabel},
	url = {http://arxiv.org/abs/2405.12762},
	doi = {10.48550/arXiv.2405.12762},
	language = {en},
	urldate = {2025-03-11},
	publisher = {arXiv},
	author = {Goulart, Paul J. and Chen, Yuwen},
	month = may,
	year = {2024},
}

@article{hofmann_performance_2023,
	title = {Performance {Assessment} of {Convex} {Low}-{Thrust} {Trajectory} {Optimization} {Methods}},
	volume = {60},
	issn = {0022-4650},
	url = {https://arc.aiaa.org/doi/10.2514/1.A35461},
	doi = {10.2514/1.A35461},
	number = {1},
	urldate = {2024-12-16},
	journal = {Journal of Spacecraft and Rockets},
	author = {Hofmann, Christian and Morelli, Andrea C. and Topputo, Francesco},
	month = jan,
	year = {2023},
	pages = {299--314},
}

@inproceedings{mao_tutorial_2018,
	title = {A {Tutorial} on {Real}-time {Convex} {Optimization} {Based} {Guidance} and {Control} for {Aerospace} {Applications}},
	url = {https://ieeexplore.ieee.org/document/8430984},
	doi = {10.23919/ACC.2018.8430984},
	urldate = {2024-04-19},
	booktitle = {2018 {Annual} {American} {Control} {Conference} ({ACC})},
	author = {Mao, Yuanqi and Szmuk, Michael and Açıkmeşe, Behçet},
	month = jun,
	year = {2018},
	pages = {2410--2416},
}

@book{nesterov_lectures_2018,
	address = {Cham},
	series = {Springer {Optimization} and {Its} {Applications}},
	title = {Lectures on {Convex} {Optimization}},
	volume = {137},
	copyright = {http://www.springer.com/tdm},
	isbn = {978-3-319-91577-7 978-3-319-91578-4},
	url = {http://link.springer.com/10.1007/978-3-319-91578-4},
	urldate = {2025-02-06},
	publisher = {Springer International Publishing},
	author = {Nesterov, Yurii},
	year = {2018},
}

@article{woffinden_navigating_2007,
	title = {Navigating the {Road} to {Autonomous} {Orbital} {Rendezvous}},
	volume = {44},
	issn = {0022-4650},
	url = {https://arc.aiaa.org/doi/10.2514/1.30734},
	doi = {10.2514/1.30734},
	number = {4},
	urldate = {2025-01-02},
	journal = {Journal of Spacecraft and Rockets},
	author = {Woffinden, David C. and Geller, David K.},
	month = jul,
	year = {2007},
	pages = {898--909},
}

@article{boone_directional_2023,
	title = {Directional {State} {Transition} {Tensors} for {Capturing} {Dominant} {Nonlinear} {Effects} in {Orbital} {Dynamics}},
	volume = {46},
	issn = {0731-5090},
	url = {https://arc.aiaa.org/doi/10.2514/1.G006910},
	doi = {10.2514/1.G006910},
	number = {3},
	urldate = {2024-02-29},
	journal = {Journal of Guidance, Control, and Dynamics},
	author = {Boone, Spencer and McMahon, Jay},
	month = mar,
	year = {2023},
	pages = {431--442},
}

@article{lu_autonomous_2013,
	title = {Autonomous {Trajectory} {Planning} for {Rendezvous} and {Proximity} {Operations} by {Conic} {Optimization}},
	volume = {36},
	issn = {0731-5090},
	url = {https://arc.aiaa.org/doi/10.2514/1.58436},
	doi = {10.2514/1.58436},
	number = {2},
	urldate = {2024-12-30},
	journal = {Journal of Guidance, Control, and Dynamics},
	author = {Lu, Ping and Liu, Xinfu},
	month = mar,
	year = {2013},
	pages = {375--389},
}

@article{malyuta_convex_2022,
	title = {Convex {Optimization} for {Trajectory} {Generation}: {A} {Tutorial} on {Generating} {Dynamically} {Feasible} {Trajectories} {Reliably} and {Efficiently}},
	volume = {42},
	issn = {1941-000X},
	shorttitle = {Convex {Optimization} for {Trajectory} {Generation}},
	url = {https://ieeexplore.ieee.org/document/9905530},
	doi = {10.1109/MCS.2022.3187542},
	number = {5},
	urldate = {2024-05-02},
	journal = {IEEE Control Systems Magazine},
	author = {Malyuta, Danylo and Reynolds, Taylor P. and Szmuk, Michael and Lew, Thomas and Bonalli, Riccardo and Pavone, Marco and Açıkmeşe, Behçet},
	month = oct,
	year = {2022},
	pages = {40--113},
}

@article{franzini_relative_2019,
	title = {Relative {Motion} {Dynamics} in the {Restricted} {Three}-{Body} {Problem}},
	volume = {56},
	issn = {0022-4650},
	url = {https://arc.aiaa.org/doi/10.2514/1.A34390},
	doi = {10.2514/1.A34390},
	number = {5},
	urldate = {2024-04-30},
	journal = {Journal of Spacecraft and Rockets},
	author = {Franzini, Giovanni and Innocenti, Mario},
	month = sep,
	year = {2019},
	pages = {1322--1337},
}

@article{valli_nonlinear_2013,
	title = {Nonlinear {Mapping} of {Uncertainties} in {Celestial} {Mechanics}},
	volume = {36},
	issn = {0731-5090},
	url = {https://arc.aiaa.org/doi/10.2514/1.58068},
	doi = {10.2514/1.58068},
	number = {1},
	urldate = {2024-02-29},
	journal = {Journal of Guidance, Control, and Dynamics},
	author = {Valli, M. and Armellin, R. and Di Lizia, P. and Lavagna, M. R.},
	month = jan,
	year = {2013},
	pages = {48--63},
}

@article{park_nonlinear_2006,
	title = {Nonlinear {Mapping} of {Gaussian} {Statistics}: {Theory} and {Applications} to {Spacecraft} {Trajectory} {Design}},
	volume = {29},
	issn = {0731-5090},
	shorttitle = {Nonlinear {Mapping} of {Gaussian} {Statistics}},
	url = {https://arc.aiaa.org/doi/10.2514/1.20177},
	doi = {10.2514/1.20177},
	number = {6},
	urldate = {2024-02-27},
	journal = {Journal of Guidance, Control, and Dynamics},
	author = {Park, Ryan S. and Scheeres, Daniel J.},
	month = nov,
	year = {2006},
	pages = {1367--1375},
}

@article{wang_minimum-fuel_2018,
	title = {Minimum-{Fuel} {Low}-{Thrust} {Transfers} for {Spacecraft}: {A} {Convex} {Approach}},
	volume = {54},
	issn = {1557-9603},
	shorttitle = {Minimum-{Fuel} {Low}-{Thrust} {Transfers} for {Spacecraft}},
	url = {https://ieeexplore.ieee.org/abstract/document/8306923},
	doi = {10.1109/TAES.2018.2812558},
	number = {5},
	urldate = {2024-12-16},
	journal = {IEEE Transactions on Aerospace and Electronic Systems},
	author = {Wang, Zhenbo and Grant, Michael J.},
	month = oct,
	year = {2018},
	pages = {2274--2290},
}

@article{acikmese_convex_2007,
	title = {Convex {Programming} {Approach} to {Powered} {Descent} {Guidance} for {Mars} {Landing}},
	volume = {30},
	issn = {0731-5090},
	url = {https://doi.org/10.2514/1.27553},
	doi = {10.2514/1.27553},
	number = {5},
	urldate = {2025-01-30},
	journal = {Journal of Guidance, Control, and Dynamics},
	author = {Acikmese, Behcet and Ploen, Scott R.},
	year = {2007},
	pages = {1353--1366},
}

@article{betts_survey_1998,
	title = {Survey of {Numerical} {Methods} for {Trajectory} {Optimization}},
	volume = {21},
	issn = {0731-5090},
	url = {https://arc.aiaa.org/doi/10.2514/2.4231},
	doi = {10.2514/2.4231},
	number = {2},
	urldate = {2025-02-06},
	journal = {Journal of Guidance, Control, and Dynamics},
	author = {Betts, John T.},
	month = mar,
	year = {1998},
	pages = {193--207},
}

@article{lu_introducing_2017,
	title = {Introducing {Computational} {Guidance} and {Control}},
	volume = {40},
	issn = {0731-5090},
	url = {https://arc.aiaa.org/doi/10.2514/1.G002745},
	doi = {10.2514/1.G002745},
	number = {2},
	urldate = {2024-12-31},
	journal = {Journal of Guidance, Control, and Dynamics},
	author = {Lu, Ping},
	month = feb,
	year = {2017},
	pages = {193--193},
}

@article{burnett_rapid_2025,
	title = {Rapid {Nonlinear} {Convex} {Guidance} {Using} a {Monomial} {Method}},
	issn = {0731-5090, 1533-3884},
	url = {https://arc.aiaa.org/doi/10.2514/1.G008512},
	doi = {10.2514/1.G008512},
	language = {en},
	urldate = {2025-02-27},
	journal = {Journal of Guidance, Control, and Dynamics},
	author = {Burnett, Ethan R. and Topputo, Francesco},
	month = feb,
	year = {2025},
	pages = {1--21},
}

@book{liberzon_calculus_2012,
	address = {Princeton ; Oxford},
	title = {Calculus of variations and optimal control theory: a concise introduction},
	isbn = {978-0-691-15187-8},
	shorttitle = {Calculus of variations and optimal control theory},
	publisher = {Princeton University Press},
	author = {Liberzon, Daniel},
	year = {2012},
	note = {OCLC: ocn772627669},
}

@article{giorgilli_methods_2011,
	title = {Methods of algebraic manipulation in perturbation theory},
	volume = {3},
	url = {https://ui.adsabs.harvard.edu/abs/2011WSAAA...3..147G},
	doi = {10.48550/arXiv.1303.7398},
	urldate = {2024-03-10},
	journal = {Workshop Series of the Asociacion Argentina de Astronomia},
	author = {Giorgilli, A. and Sansottera, M.},
	month = jan,
	year = {2011},
	pages = {147--183},
}

@techreport{lee_white_2019,
	title = {White {Paper}: {Gateway} {Destination} {Orbit} {Model}: {A} {Continuous} 15 {Year} {NRHO} {Reference} {Trajectory}},
	shorttitle = {White {Paper}},
	url = {https://ntrs.nasa.gov/citations/20190030294},
	urldate = {2024-05-16},
	institution = {NASA Technical Management},
	author = {Lee, David E.},
	month = aug,
	year = {2019},
}

@article{spreen_automated_2017,
	title = {Automated {Node} {Placement} {Capability} for {Spacecraft} {Trajectory} {Targeting} {Using} {Higher}-{Order} {State} {Transition} {Matrices}},
	language = {en},
	journal = {AAS/AIAA Astrodynamics Specialist Conference},
	author = {Spreen, Christopher and Howell, Kathleen},
	year = {2017},
	pages = {21--24},
}

@article{gimeno_numerical_2023,
	title = {Numerical integration of high-order variational equations of {ODEs}},
	volume = {442},
	issn = {00963003},
	url = {https://linkinghub.elsevier.com/retrieve/pii/S0096300322008116},
	doi = {10.1016/j.amc.2022.127743},
	language = {en},
	urldate = {2024-12-30},
	journal = {Applied Mathematics and Computation},
	author = {Gimeno, Joan and Jorba, Angel and Jorba-Cuscó, Marc and Miguel, Narcís and Zou, Maorong},
	month = apr,
	year = {2023},
	pages = {127743},
}

@article{liu_survey_2017,
	title = {Survey of convex optimization for aerospace applications},
	volume = {1},
	issn = {2522-0098},
	url = {https://doi.org/10.1007/s42064-017-0003-8},
	doi = {10.1007/s42064-017-0003-8},
	language = {en},
	number = {1},
	urldate = {2024-12-16},
	journal = {Astrodynamics},
	author = {Liu, Xinfu and Lu, Ping and Pan, Binfeng},
	month = sep,
	year = {2017},
	pages = {23--40},
}

@incollection{hofmann_performance_2021,
	series = {{AIAA} {SciTech} {Forum}},
	title = {On the {Performance} of {Discretization} and {Trust}-{Region} {Methods} for {On}-{Board} {Convex} {Low}-{Thrust} {Trajectory} {Optimization}},
	url = {https://arc.aiaa.org/doi/10.2514/6.2022-1892},
	urldate = {2024-12-16},
	booktitle = {{AIAA} {SCITECH} 2022 {Forum}},
	publisher = {American Institute of Aeronautics and Astronautics},
	author = {Hofmann, Christian and Morelli, Andrea C. and Topputo, Francesco},
	month = dec,
	year = {2021},
	doi = {10.2514/6.2022-1892},
}

@book{boyd_convex_2004,
	address = {Cambridge, UK ; New York},
	title = {Convex optimization},
	isbn = {978-0-521-83378-3},
	language = {en},
	publisher = {Cambridge University Press},
	author = {Boyd, Stephen P. and Vandenberghe, Lieven},
	year = {2004},
}

@article{lizy-destrez_rendezvous_2019,
	title = {Rendezvous {Strategies} in the {Vicinity} of {Earth}-{Moon} {Lagrangian} {Points}},
	volume = {5},
	issn = {2296-987X},
	url = {https://www.frontiersin.org/article/10.3389/fspas.2018.00045/full},
	doi = {10.3389/fspas.2018.00045},
	language = {en},
	urldate = {2024-04-29},
	journal = {Frontiers in Astronomy and Space Sciences},
	author = {Lizy-Destrez, Stephanie and Beauregard, Laurent and Blazquez, Emmanuel and Campolo, Antonino and Manglativi, Sara and Quet, Victor},
	month = jan,
	year = {2019},
	pages = {45},
}

@article{scorsoglio_relative_2023,
	title = {Relative motion guidance for near-rectilinear lunar orbits with path constraints via actor-critic reinforcement learning},
	volume = {71},
	issn = {02731177},
	url = {https://linkinghub.elsevier.com/retrieve/pii/S0273117722007207},
	doi = {10.1016/j.asr.2022.08.002},
	language = {en},
	number = {1},
	urldate = {2024-04-29},
	journal = {Advances in Space Research},
	author = {Scorsoglio, Andrea and Furfaro, Roberto and Linares, Richard and Massari, Mauro},
	month = jan,
	year = {2023},
	pages = {316--335},
}

@article{short_stretching_2015,
	title = {Stretching in phase space and applications in general nonautonomous multi-body problems},
	volume = {122},
	issn = {1572-9478},
	url = {https://doi.org/10.1007/s10569-015-9617-4},
	doi = {10.1007/s10569-015-9617-4},
	language = {en},
	number = {3},
	urldate = {2024-04-10},
	journal = {Celestial Mechanics and Dynamical Astronomy},
	author = {Short, Cody R. and Blazevski, Daniel and Howell, Kathleen C. and Haller, George},
	month = jul,
	year = {2015},
	pages = {213--238},
}

@article{mao_successive_2017,
	series = {20th {IFAC} {World} {Congress}},
	title = {Successive {Convexification} of {Non}-{Convex} {Optimal} {Control} {Problems} with {State} {Constraints}},
	volume = {50},
	issn = {2405-8963},
	url = {https://www.sciencedirect.com/science/article/pii/S2405896317312405},
	doi = {10.1016/j.ifacol.2017.08.789},
	number = {1},
	urldate = {2024-04-05},
	journal = {IFAC-PapersOnLine},
	author = {Mao, Yuanqi and Dueri, Daniel and Szmuk, Michael and Açıkmeşe, Behçet},
	month = jul,
	year = {2017},
	pages = {4063--4069},
}

@article{park_nonlinear_2007,
	title = {Nonlinear {Semi}-{Analytic} {Methods} for {Trajectory} {Estimation}},
	volume = {30},
	doi = {10.2514/1.29106},
	journal = {Journal of Guidance Control and Dynamics - J GUID CONTROL DYNAM},
	author = {Park, Ryan and Scheeres, D.},
	month = nov,
	year = {2007},
	pages = {1668--1676},
}

@article{survey-topputo,
author = {Topputo, F. and Zhang, C.},
title = {Survey of Direct Transcription for Low-Thrust Space Trajectory Optimization with Applications},
journal = {Abstract and Applied Analysis},
volume = {2014},
number = {1},
pages = {851720},
doi = {https://doi.org/10.1155/2014/851720},
url = {https://onlinelibrary.wiley.com/doi/abs/10.1155/2014/851720},
eprint = {https://onlinelibrary.wiley.com/doi/pdf/10.1155/2014/851720},
year = {2014}
}

@article{impulsiveSTMSCP_rizza,
author = {Antonio Rizza and Francesco Topputo and Simone D'Amico},
title = {Goal-oriented asteroid mapping under uncertainties using Sequential Convex Programming},
booktitle = {AIAA SCITECH 2024 Forum},
pages={1990},
year={2024},
doi = {10.2514/6.2024-1990},
URL = {https://arc.aiaa.org/doi/abs/10.2514/6.2024-1990},
eprint = {https://arc.aiaa.org/doi/pdf/10.2514/6.2024-1990},
}

@misc{benedikter2019convexoptimizationlinearimpulsive,
      title={Convex Optimization of Linear Impulsive Rendezvous}, 
      author={Boris Benedikter and Alessandro Zavoli},
      year={2019},
      eprint={1912.08038},
      archivePrefix={arXiv},
      primaryClass={math.OC},
      url={https://arxiv.org/abs/1912.08038}, 
}

@book{primer-vector-1963,
author="Lawden, Derek F",
title="Optimal trajectories for space navigation",
publisher="Butterworths",
year="1963",
series="Butterworths mathematical texts",
}

@article{directOptmization,
author = {Hargraves, C.R. and Paris, S.W.},
title = {Direct trajectory optimization using nonlinear programming and collocation},
journal = {Journal of Guidance, Control, and Dynamics},
volume = {10},
number = {4},
pages = {338-342},
year = {1987},
doi = {10.2514/3.20223},

URL = { 
    
        https://doi.org/10.2514/3.20223
    
    

},
eprint = { 
    
        https://doi.org/10.2514/3.20223
    
    

}

}

@article{wittig_propagation_2015,
	title = {Propagation of large uncertainty sets in orbital dynamics by automatic domain splitting},
	volume = {122},
	issn = {1572-9478},
	url = {https://doi.org/10.1007/s10569-015-9618-3},
	doi = {10.1007/s10569-015-9618-3},
	language = {en},
	number = {3},
	urldate = {2024-03-24},
	journal = {Celestial Mechanics and Dynamical Astronomy},
	author = {Wittig, Alexander and Di Lizia, Pierluigi and Armellin, Roberto and Makino, Kyoko and Bernelli-Zazzera, Franco and Berz, Martin},
	month = jul,
	year = {2015},
	pages = {239--261},
}

@inproceedings{boone_optimal_2021,
	title = {Optimal {Maneuver} {Targeting} {Using} {State} {Transition} {Tensors} with {Variable} {Time}-of-{Flight}},
    author={Boone, Spencer and McMahon, Jay},
    booktitle={31st AAS/AIAA Space Flight Mechanics Meeting},
    month = feb,
    year={2021}
}

@article{di_lizia_high_2014,
	title = {High order optimal control of space trajectories with uncertain boundary conditions},
	volume = {93},
	issn = {0094-5765},
	url = {https://www.sciencedirect.com/science/article/pii/S0094576513002397},
	doi = {10.1016/j.actaastro.2013.07.007},
	urldate = {2024-03-08},
	journal = {Acta Astronautica},
	author = {Di Lizia, P. and Armellin, R. and Bernelli-Zazzera, F. and Berz, M.},
	month = jan,
	year = {2014},
	pages = {217--229},
}

@article{kolda_tensor_2009,
	title = {Tensor {Decompositions} and {Applications}},
	volume = {51},
	issn = {0036-1445, 1095-7200},
	url = {http://epubs.siam.org/doi/10.1137/07070111X},
	doi = {10.1137/07070111X},
	language = {en},
	number = {3},
	urldate = {2024-03-01},
	journal = {SIAM Review},
	author = {Kolda, Tamara G. and Bader, Brett W.},
	month = aug,
	year = {2009},
	pages = {455--500},
}

\end{document}